\documentclass[12pt,reqno,a4paper]{article}
\usepackage{lineno,hyperref}
\usepackage{tikz}
\usetikzlibrary{shapes}
\usepackage{amsfonts,amsmath,times,amssymb}
\usepackage{cool}
\usepackage{xparse}
\usepackage{verbatim}
\usepackage{multirow}

\newtheorem{theorem}{Theorem}
\newtheorem{lemma}{Lemma}
\newtheorem{prop}{Proposition}

\newcommand{\clique}{\mathcal{C}}
\newcommand{\apol}{\mathcal{A}}
\newcommand{\matA}{\boldsymbol{A}}
\newcommand{\Given}{\, \Bigg{\vert} \,}
\renewcommand{\E}{\mathbb{E}}
\newcommand{\indicator}{{\bf 1}}
\newcommand{\convas}{\, \overset{a.s.}{\longrightarrow} \,}
\newcommand{\Prob}{\mathbb{P}}

\newcommand{\F}{\mathcal{F}}
\newcommand{\polya}{P\'{o}lya}
\NewDocumentCommand{\ktree}{mg}{\ensuremath{\mathcal{T}^{(#1)}_{#2}}}
\NewDocumentCommand{\ext}{mg}{\ensuremath{\mathcal{E}^{(#1)}_{#2}}}
\NewDocumentCommand{\dsq}{mg}{\ensuremath{\mathcal{D}^{(#1)}_{#2}}}

\begin{document}
\begin{center}
	{\Large \bf Characterizing several properties of high-dimensional random Apollonian networks}
	
	\bigskip
	{\large \bf Panpan Zhang\footnote{Department of Biostatistics, Epidemiology and Informatics, Perelman School of Medicine, University of Pennsylvania, Philadelphia, PA 19104, U.S.A.}}
	
	\bigskip
	\bf{\today}
\end{center}

\bigskip\noindent
{\bf Abstract.} In this article, we investigate several properties 
of high-dimensional random Apollonian networks (HDRANs), including 
two types of degree profiles, the small-world effect (clustering 
property), sparsity, and three distance-based metrics. The 
characterizations of degree profiles are based on several rigorous 
mathematical and probabilistic methods, such as a two-dimensional 
mathematical induction, analytic combinatorics, and \polya\ urns, 
etc. The small-world property is uncovered by a well-developed 
measure---local clustering coefficient, and the sparsity is assessed
by a proposed Gini index. Finally, we look into three distance-based 
properties; they are total depth, diameter and Wiener index.


Secondary: 05C07; 05C12
  

\section{Introduction}
Due to the surge of interest in networks, such as the Internet 
(e.g., The Internet Mapping Project by Hal Burch and Bill Cheswick), 
resistor networks~\cite{newman2004finding}, the World Wide Web 
(WWW)~\cite{barabasi1999emergence}, and social networks (e.g., 
friendship 
network~\cite{moody2001race}), a plethora of network models have 
been 
proposed and studied in the last several decades. In this paper, we 
investigate a network model that recently caught researchers' 
attention---{\em Apollonian networks} (ANs). ANs arise from the 
problem of space-filling packing of spheres, proposed by the ancient 
Greek mathematician Apollonius of Perga. ANs possess a variety of 
typical network characteristics, 
which are summarized in the title of~\cite{andrade2005apollonian}: 
scale free, 
small world, Euclidean, space filling and matching graphs. Each of 
these phrases is a significant area of modern network research 
itself. In practice, ANs have found a lot of applications in 
different scientific disciplines~\cite{almeida2013quantum, 
	huang2006walks, 
	lima2012nonequilibrium, pawela2015generalized, serva2013ising, 
	silva2013critical, souza2013discrete, wong2010partially, 
	xu2008coherent}. The wide application of this class of networks 
motivates us to conduct the present research.

The counterpart of AN in the field of random network analysis is 
called {\em Random Apollonian Network} (RAN). The study of RAN first 
appeared in~\cite{zhou2005maximal}, where the power-law and the 
clustering 
coefficient were investigated. Since then, many more properties of 
RANs have been uncovered by applied mathematicians and probablists: 
The degree distribution was characterized by~\cite{frieze2014some}; 
the diameter was calculated by~\cite{ebrahimzadeh2014onlongest, 
	frieze2014some}; the length 
of the longest path in RANs was determined 
by~\cite{collevecchio2016longest, cooper2015long, 
	ebrahimzadeh2014onlongest}. All these research papers, however, 
	only 
focused on planar RANs, the evolution of which is based on 
continuing triangulation. Triangulated RANs are a special class of 
(more general) RANs with network {\em index} taking value $3$. It 
can be 
shown that triangulated RANs are maximal planar graphs by the {\em 
	Kuratowski 	criterion}~\cite{kuratowski1930remarques}. It is 
	evident 
that there is an underlying theory of {\em preferential attachment} 
(PA)~\cite{massen2007preferential} in the evolution of RANs, where 
PA is a critical manifestation in social sciences, rendering the 
potential application of RANs in a wider range of fields (than what 
has been found in the literature).

In this paper, we consider 
a class of high-dimensional networks generalized from triangulated 
RANs, i.e., high-dimensional random Apollonian Networks (HDRANs) 
that refer to the RANs with a general network index $k \ge 3$. 
HDRANs were first introduced by~\cite{zhang2006high}, where an 
iterative 
algorithm was designed to characterize several network properties 
including degree distribution, clustering coefficient and diameter. 
The exact 
degree distribution of a vertex with a fixed label and the total 
weight (a macro metric) were determined 
by~\cite{zhang2016thedegree}. A follow-up study embedding RANs into 
continuous time was given by~\cite{zhang2016distributions}. To the 
best of our knowledge, there is almost no other work has been done 
for 
HDRANs in the literature. 

The goal of this paper is to give a comprehensive study of HDRANs, 
with specific focus on the investigation of several network 
properties of common interest by utilizing some well-developed 
methods; for instance, stochastic recurrences and \polya\ urns. 
Some of the results, such as the degree distribution, are 
directly extended from their counterparts for triangulated RANs. For 
better readability, only the results and the underlying theory are 
presented in the main body of the paper, but the associated 
mathematical derivations are given in the appendix. A couple of 
novel network properties, such as the sparsity and the 
depth of HDRANs, are rigorously uncovered as well. Details will be 
given in the 
sequel. 

The rest of the paper is 
organized as follows. In Section~\ref{Sec:evolution}, we briefly 
review the evolutionary process of HDRANs as well as some basic 
graph invariants thereof. In the next five sections, the main 
results of the analyzed properties are given; see a summary in 
Table~\ref{Table:summary}.
\begin{table}[h!]
	\begin{center}
		\renewcommand*{\arraystretch}{1.2}
		\begin{tabular}{|c|c|c|}
			\hline
			Section & Property & Method(s)
			\\ \hline
			$3$ & Degree profile {\rm I} & Two-dimensional induction 
			(extended from~\cite{frieze2014some})
			\\ \hline
			\multirow{2}{*}{$4$} & \multirow{2}{*}{Degree profile 
				{\rm II}} & Analytic 
			combinatorics~\cite{flajolet2006some}
			\\ & & Triangular urns~\cite{zhang2016thedegree}
			\\ \hline
			$5$ & Small world & Local clustering coefficient
			\\ \hline
			$6$ & Sparsity & A proposed Gini index
			\\ \hline
			\multirow{3}{*}{$7$} & Total depth & Recurrence methods
			\\ & Diameter & Results directly 
			from~\cite{cooper2014theheight}
			\\ & The Wiener index & Numeric experiments
			\\ \hline
		\end{tabular}
		\caption{Summary of the article}
		\label{Table:summary}
	\end{center}
\end{table}
In Section~\ref{Sec:concluding}, we give some concluding remarks and 
propose some future work.

\section{Evolution of random Apollonian networks}
\label{Sec:evolution}

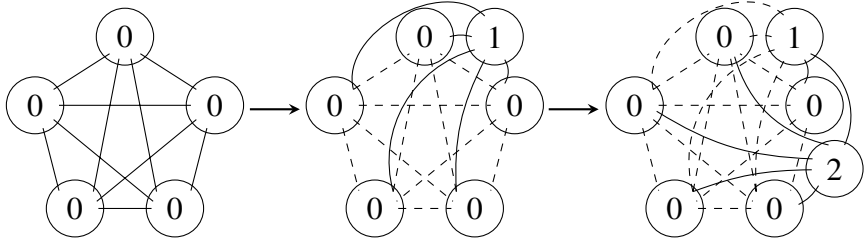
\begin{figure}[tbh]
	\begin{center}
		\begin{tikzpicture}[scale=2.63]
			\draw
			(-0.25,0) node [circle=0.1,draw] {0}
			(0.25,0) node [circle=0.1,draw] {0}
			(0,0.866) node [circle=0.1,draw] {0}
			(-0.45, 0.52) node [circle=0.1,draw] {0}
			(0.45, 0.52) node [circle=0.1,draw] {0}
			(-0.133,0) -- (0.133,0)
			(-0.33,0.52) -- (0.33,0.52)
			(-0.158, 0.085) -- (-0.03, 0.755)
			(0.158, 0.085) -- (0.03, 0.755)
			(-0.41, 0.405) -- (-0.33, 0.09)
			(0.41, 0.405) -- (0.33, 0.09)
			(-0.357, 0.448) -- (0.141, 0.036)
			(0.357, 0.448) -- (-0.141, 0.036)
			(-0.357, 0.6) -- (-0.077, 0.78)
			(0.357, 0.6) -- (0.077, 0.78);
			\draw [->,>=stealth,thick] (0.625,.5) -- +(0.25,0);
			\draw
			(1.25,0) node [circle=0.1,draw] {0}
			(1.75,0) node [circle=0.1,draw] {0}
			(1.5,0.866) node [circle=0.1,draw] {0}
			(1.05, 0.52) node [circle=0.1,draw] {0}
			(1.95, 0.52) node [circle=0.1,draw] {0}
			(1.85, 0.866) node [circle=0.1,draw] {1};
			\draw[dashed]
			(1.367,0) -- (1.633,0)
			(1.17,0.52) -- (1.83,0.52)
			(1.342, 0.085) -- (1.47, 0.755)
			(1.658, 0.085) -- (1.53, 0.755)
			(1.09, 0.405) -- (1.17, 0.09)
			(1.91, 0.405) -- (1.83, 0.09)
			(1.143, 0.448) -- (1.641, 0.036)
			(1.857, 0.448) -- (1.359, 0.036)
			(1.143, 0.6) -- (1.423, 0.78)
			(1.857, 0.6) -- (1.577, 0.78);
			\draw
			(1.625, 0.866) to[bend left=15] (1.735, 0.866)
			(1.143, 0.6) to[bend left=75] (1.8, 0.97)
			(1.342, 0.085) to[bend left=45] (1.75, 0.8)
			(1.658, 0.085) to[bend left=15] (1.8, 0.75)
			(1.9, 0.63) to[bend right=30] (1.9, 0.75);			
			\draw [->,>=stealth,thick] (2.12,.5) -- +(0.25,0);
			\draw
			(2.75,0) node [circle=0.1,draw] {0}
			(3.25,0) node [circle=0.1,draw] {0}
			(3,0.866) node [circle=0.1,draw] {0}
			(2.55, 0.52) node [circle=0.1,draw] {0}
			(3.45, 0.52) node [circle=0.1,draw] {0}
			(3.35, 0.866) node [circle=0.1,draw] {1}
			(3.55, 0.2) node [circle=0.1,draw] {2};
			\draw[dashed]
			(2.867,0) -- (3.133,0)
			(2.67,0.52) -- (3.33,0.52)
			(2.842,0.085) -- (2.97, 0.755)
			(3.158, 0.085) -- (3.03, 0.755)
			(2.59, 0.405) -- (2.67, 0.09)
			(3.41, 0.405) -- (3.33, 0.09)
			(2.643, 0.448) -- (3.141, 0.036)
			(3.357, 0.448) -- (2.859, 0.036)
			(2.643, 0.6) -- (2.923, 0.78)
			(3.357, 0.6) -- (3.077, 0.78)
			(3.125, 0.866) to[bend left=15] (3.235, 0.866)
			(2.643, 0.6) to[bend left=75] (3.3, 0.97)
			(2.842, 0.085) to[bend left=45] (3.25, 0.8)
			(3.158, 0.085) to[bend left=15] (3.3, 0.75);
			\draw
			(3.4, 0.63) to[bend right=30] (3.4, 0.75)
			(3.6, 0.32) to [bend right=45] (3.43, 0.78)
			(3.52, 0.32) to [bend left=30] (3.05, 0.765)
			(3.44, 0.25) to [bend left=15] (2.66,0.49)
			(3.436, 0.195) to [bend right=10] (2.842, 0.085)
			(3.47, 0.115) to [bend left=15] (3.37, 0.02);	
		\end{tikzpicture}
		\caption{An example of the evolution of a HDRAN of index 5 
			in two steps; active cliques are those containing at 
			least 
			one solid edge.}
		\label{Fig:evol}
	\end{center} 
\end{figure}

In this section, we review the evolution of a RAN of index $k \ge 
3$. At time $n = 0$, we start with a {\em complete 
	graph}\footnote{In graph theory, a complete graph is a graph 
	such 
	that each pair of vertices therein is connected by an edge. A 
	complete graph on $k$ vertices is also called $k$-clique or 
	$k$-simplex. We shall interchangeably use these terms through 
	the 
	manuscript.} on $k$ vertices all of which are labeled with $0$. 
	At 
each subsequent time point $n \ge 1$, a $k$-clique is chosen 
uniformly at random among all active cliques in the network. A new 
vertex labeled with $n$ is linked by $k$ edges to all the vertices 
of the chosen clique. Then, the recruiting clique is deactivated. An 
explanatory example of a RAN with index $k = 5$ is given in 
Figure~\ref{Fig:evol}.

According to the evolutionary process described above, we obtain 
some basic and deterministic graph invariants of a RAN with index 
$k$ at time $n$: the number of vertices $V_{n}^{(k)} = k + n$, the 
number of edges $E_{n}^{(k)} = k + nk$, and the number of active 
cliques $\clique_n^{(k)} = 1 + (k - 1)n$. We note that RANs of 
indices $1$ and $2$ are not considered in this paper, as their 
structure lacks research interest. A RAN of index $1$ at time $n$ is 
a single vertex labeled with $n$, while a RAN of index $2$ at time 
$n$ is a path of length $n$.

\section{Degree profile {\rm I}}
\label{Sec:degree}
In this section, we investigate the degree profile of a RAN of index 
$k \ge 3$. The random variable of prime interest is $X_{n, 
	j}^{(k)}$, the number of vertices of degree $j$ in a RAN of 
	index 
$k$ at time $n$, for $j \ge k$, where the boundary condition arises 
from the natural lower bound of the degree of vertices in 
RANs\footnote{Upon joining into the network, every newcomer is 
	connected with $k$ existing vertices, leading to minimal 
	possible 
	degree $k$.}. It is also worthy of noting that the natural upper 
bound for $j$ at time $n$ is $k + n - 1$. 

The degree random variable that we consider in this section is 
different from that investigated in~\cite{zhang2016thedegree}, and 
the methods 
developed in~\cite{zhang2016thedegree} are not amenable to this 
study, which will 
be explained in detail in the sequel. To distinguish the two kinds 
of degree profiles, we call the one discussed in this section degree 
profile {\rm I}. Specifically, we present two results of $X_{n, 
	j}^{(k)}$, which are respectively shown in 
Theorems~\ref{Thm:L1bound} and~\ref{Thm:pbound}. In 
Theorem~\ref{Thm:L1bound}, we prove that the difference between the 
expectation of $X_{n, j}^{(k)}$ and a linear function of $n$ is 
uniformly bounded, where the bound is determined. In 
Theorem~\ref{Thm:pbound}, we show that $X_{n, j}^{(k)}$ concentrates 
on its expectation with high probability, i.e., a focusing property. 

\begin{theorem}
	\label{Thm:L1bound}
	Let $X_{n, j}^{(k)}$ be the number of vertices of degree $j$ in 
	a RAN of index $k$ at time $n$, for $j \ge k$. For each $n \in 
	\mathbb{N}$ and any $k \ge 3$, there exists a constant $b_{j, 
		k}$ such that 
	\begin{equation}
		\label{Eq:degreeL1}
		\left|\E \left[X_{n, j}^{(k)}\right] - b_{j, k} \, n\right| 
		\le \frac{2k^2}{2k - 1}.
	\end{equation}
	In particular, we have $b_{j, k} = \frac{\Gamma(j)\Gamma(2k - 
		1)}{\Gamma(j + k) \Gamma(k - 1)}$.
\end{theorem} 

The proof of Theorem~\ref{Thm:L1bound} is based on an elementary 
mathematical tool---induction. As suggested 
in~\cite{frieze2014some}, we 
split the cases of $j = k$ and $j > k$ in the proof. For the case of 
$j = k$, we apply the traditional mathematical induction directly, 
whereas we develop a two-dimensional induction based on an infinite 
triangular array for the case of $j > k$. For the better readability 
of the paper, we present the major steps of the proof in 
Appendix~\ref{App:L1bound}. 

In the proof of Theorem~\ref{Thm:L1bound}, we show that the mean of 
$X_{n, j}^{(k)}$ scaled by $n$ converges to $b_{j, k}$ when $n$ is 
large. When $j$ goes to infinity as well, we discover that $b_{j, k} 
\sim j^{-k}$ according to the {\em Stirling's approximation}. This 
implies that the degree distribution in HDRANs follows a {\em 
	power-law} property, where the exponent is the network index 
	$k$. 
Consequently, HDRANs are {\em scale-free} networks. The power-law 
property for planar RANs (i.e., $k = 3$) has been recovered 
in~\cite{zhou2005maximal} numerically and in~\cite{frieze2014some} 
analytically.

In addition, we are interested in the deviation of the random 
variable $X_{n, j}$ from its expectation. In 
Thoerem~\ref{Thm:pbound}, we develop a Chebyshev-type inequality. 

\begin{theorem}
	\label{Thm:pbound}
	Let $X_{n, j}^{(k)}$ be the number of vertices of degree $j$ in 
	a RAN of index $k$ at time $n$, for $j \ge k$. For any $\lambda> 
	0$, we have
	$$\Prob\left(\left|X_{n, j}^{(k)} - \E\left[X_{n, j}^{(k)} 
	\right]\right| \ge \lambda \right) \le e^{-\lambda^2/(8kn)}.$$
\end{theorem}

The proof of Theorem~\ref{Thm:pbound} is presented in 
Appendix~\ref{App:pbound}. The main idea is to employ the {\em 
	Azuma-Hoeffding inequality}~\cite{azuma1967weighted} based on a 
martingale 
sequence. We remark that the exact same concentration result is 
found for {\em random $k$-trees}~\cite{gao2009thedegree}. The author 
of~\cite{gao2009thedegree} tackled the problem by using the methods 
from tree 
realization theory. The intrinsic reason of the identicality is 
similarity in the evolutionary processes of HDRANs with index 
$k$ and random $k$-trees.

Before ending this section, we would like to point out that the 
methods in the proofs of Theorems~\ref{Thm:L1bound} 
and~\ref{Thm:pbound} are extended from the ideas 
in~\cite{frieze2014some}. 
The results for planar RANs (a special case for $k = 3$) can be 
found in~\cite[Theorem 1.1]{frieze2014some}.

\section{Degree profile {\rm II}}
\label{Sec:degree2}
Another type of degree profile that we look into is node-specified. 
Let $D_{n, j}^{(k)}$ denote the degree of the node labeled with $j$ 
in a HDRAN of index $k$ at time $n$. This property was investigated 
in~\cite{zhang2016thedegree}, where the growth of HDRANs was 
represented by a 
two-color \polya\ urn scheme~\cite{mahmoud2009polya}. \polya\ urn 
appears to 
be an appropriate model since it successfully captures the 
evolutionary characteristics of highly dependent structures.

Noticing that the degree of a vertex is equal to the number of 
cliques incident with it, the authors of~\cite{zhang2016thedegree} 
introduced a 
color code such that the active cliques incident with the node 
labeled with $j$ were colored white, while all the rest were colored 
blue. The associated urn scheme is governed by the {\em replacement 
	matrix}
$$\begin{pmatrix}
	k - 2 & 1
	\\ 0 & k - 1
\end{pmatrix}.$$  
This replacement matrix is triangular, so the associated \polya\ urn 
is called {\em triangular urn}. This class of urns has been 
extensively studied in~\cite{flajolet2006some, janson2006limit, 
	zhang2015explicit}. The next 
proposition specifies the exact distribution of $D_{j, n}^{(k)}$ as 
well as its moments.
\begin{prop}
	\label{Thm:degreedist}
	Let $D_{n, j}^{(k)}$ be the degree of the node labeled with $j$ 
	in a RAN of index $k$ at time $n$, for $n \ge j$. The 
	distribution of $D_{n, j}^{(k)}$ is given by
	\begin{align*}
		\Prob\left(D_{n, j}^{(k)} = k + \delta\right) &= 
		\frac{\Gamma(n - j + 1)\Gamma\left(j + \frac{1}{k - 
				1}\right)}{\Gamma\left(n + \frac{1}{k - 
				1}\right)}{{\delta + 
				\frac{2}{k - 2}} \choose \delta} 
		\\ &\qquad{}\times \sum_{r = 1}^{\delta} (-1)^{\delta} 
		{\delta \choose i} {{n - 2 - \frac{k - 2}{k - 1} r} \choose 
			{n - j}},
	\end{align*}
	for $\delta = 1, 2, \ldots, n - j$. The $s$-th moment of $D_{n, 
		j}^{(k)}$ is
	\begin{align}
		\E\left[\left(D_{n, j}^{(k)}\right)^{s\,}\right] &= 
		\frac{1}{(k - 2)^s} \left[(k(k - 3))^s + \sum_{r = 1}^{s}{s 
			\choose r} \frac{(k(k - 3))^{s - r}(k - 
			2)^r}{\left\langle j 
			+ 1/(k - 1)\right\rangle_{n - j}}\right. \nonumber
		\\ &\qquad{}\left.\times \sum_{i = 1}^{r} (-1)^{r - i} {r 
			\brace i} \left\langle \frac{k}{k - 2} \right\rangle_i 
		\left\langle j + \frac{1}{k - 1} + \frac{k - 2}{k - 1}i 
		\right\rangle_{n - j}\right], \label{Eq:degreemoment}
	\end{align}
	where $\langle \cdot \rangle_{\cdot}$ represents the Pochhammer 
	symbol of rising factorial, and ${{\cdot} \brace {\cdot}}$ 
	represents Stirling numbers of the second kind.
\end{prop}

The probability distribution function of $D_{n, j}^{(k)}$ is 
obtained by exploiting the results in~\cite[Proposition 
14]{flajolet2006some}, and the moments are recovered 
from~\cite[Proposition 
1]{zhang2016thedegree}. The asymptotic moments of $D_{n, j}^{(k)}$ 
are obtained 
directly by applying the Stirling's approximation to 
Equation~(\ref{Eq:degreemoment}); namely,

$$\E\left[\left(\frac{D_{n, j}^{(k)}}{n^{(k - 2)/(k - 1)}}\right)^{s 
	\,}\right] = \frac{\Gamma\left(j + \frac{1}{k - 1}\right) 
	\Gamma\left( s + \frac{k}{k - 2}\right)}{\Gamma\left(j + \frac{k 
	- 
		2}{k - 1}s + \frac{1}{k - 1}\right) \Gamma\left(\frac{k}{k - 
		2}\right)}.$$

In particular, the asymptotic mean of $D_{n, j}^{(k)}$ is given by

$$\E\left[D_{n, j}^{(k)}\right] \sim \frac{\frac{k}{k - 2} 
	\Gamma\left(j + \frac{1}{k - 1}\right)}{\Gamma(j + 1)} \, n^{(k 
	- 
	2)/(k - 1)},$$
implying a phase transition in $j = j(n)$:
$$\E\left[D_{n, j}^{(k)}\right] \sim
\begin{cases}
	\frac{k}{k - 2} \left(\frac{n}{j(n)}\right)^{(k - 2)/(k - 1)}, 
	\qquad &j = o(n),
	\\ \frac{k}{k - 2} \left(k - 3 + \alpha^{-(k - 2)/(k - 1)} 
	\right), \qquad &j \sim \alpha n, 
\end{cases}
$$
for some $\alpha > 0$.

\section{Small-world}
\label{Sec:smallworld}
In 
this section, we look into the {\em small-world} property of HDRANs. 
The term ``small-world'' was first coined 
by~\cite{watts1998collective}. In the 
paper, the authors suggested to use the average of {\em local 
	clustering coefficients} to assess the small-world effect of a 
network; that is,
$$\hat{C}(n) = \frac{1}{V} \sum_{v} C_v(n),$$
where $V = V_n^{(k)}$ denotes the total number of vertices and 
$C_v(n)$ is the 
local clustering coefficient of vertex~$v$ at time $n$. The local 
clustering coefficient of vertex $v$ is defined as the proportion of 
the number of edges in the {\em open neighborhood} of $v$, i.e.,
$$C_v(n) = \frac{|\{e_{uw}, u, w \in 
	\mathcal{N}_v(n)\}|}{|\mathcal{N}_v(n)|(|\mathcal{N}_v(n) - 
	1|)/2},$$
where $\mathcal{N}_v(n)$ is the open neighborhood of $v$ at time 
$n$, $e_{ij}$ denotes an edge between vertices $i$ and $j$, and 
$|\cdot|$ represents the cardinality of a set.

For each newcomer, say $v^{*}$, to an HDRAN of index $k$, the open 
neighborhood of $v^{*}$ is comprised by $k$ vertices of the simplex 
chosen for recruiting $v^*$. Thus, the order of the open 
neighborhood of $v^{*}$ is $k$, and the number of edges in the 
neighborhood of $v^{*}$ is ${k \choose 2}$. Upon the first 
appearance of $v^{*}$ in the network, the degree of $v^{*}$, denoted 
$d_{v^{*}}(n)$, is $k$. As an active simplex containing $v^{*}$ is 
selected for recruiting a newcomer in any subsequent time point, 
$d_{v^{*}}(n)$ increases by $1$, and the number of edges of the 
neighborhood of $v^{*}$ increases by $k - 1$. In general, for a 
vertex $v$ of ${\rm deg}_v(n) = j$ at time $n$, the clustering 
coefficient is given by
\begin{equation*}
	C_v(n) = \frac{(k - 1)(j - k) + {k \choose 2}}{{j \choose 2}} = 
	\frac{(k - 1)(2j - k)}{j(j - 1)}.
\end{equation*}
Accordingly, the clustering coefficient of the entire network at 
time $n$ is
\begin{equation*}
	\hat{C}(n) = \frac{1}{n + k} \sum_{v}C_v(n) = \sum_{j = k}^{k + 
		n - 1} \frac{(k - 1)(2j - k)}{j(j - 1)} \times  \frac{X_{n, 
			j}^{(k)}}{n + k},
\end{equation*}
where $X_{n, j}^{(k)}$ denotes the number of vertices of degree $j$ 
in the network at time $n$. When the network is large (i.e., $n \to 
\infty$), the asymptotic clustering coefficient is given by
\begin{equation*}
	\hat{C}(\infty) \approx \sum_{j = k}^{\infty} \frac{(k - 1)(2j 
		- k)}{j(j - 1)} \lim_{n \to \infty} \frac{\E \left[X_{n, 
			j}^{(k)}\right]}{n + k} = \sum_{j = k}^{\infty} \frac{(k 
			- 1)(2j 
		- k)}{j(j - 1)} 
	\frac{\Gamma(j) \Gamma(2k - 1)}{\Gamma(j + k) \Gamma(k - 1)},
\end{equation*}
where the second equality in the last display holds according to 
Theorem~\ref{Thm:L1bound}. We simplify the expression of 
$\hat{C}(\infty)$ by applying several algebraic results of gamma 
function, and get
\begin{equation*}
	\label{Eq:asymcc}
	\hat{C}(\infty) \approx \frac{(k - 1)\Gamma(2k - 1)}{\Gamma(k - 
		1)} \sum_{j = k}^{\infty} \frac{(2j - k) \Gamma(j - 1)}{j \, 
		\Gamma(j + k)} = \frac{(k - 1)\Gamma(2k - 1)}{\Gamma(k - 1)} 
	\sum_{j = 
		k}^{\infty} \left(\frac{2 \, \Gamma(j - 1)}{\Gamma(j + k)} - 
	\frac{k \, \Gamma(j - 1)}{j \, \Gamma(j + k))} \right).
\end{equation*}
We evaluate the two terms in the summand one after another. The 
first sum is given by
$$\sum_{j = k}^{\infty} \frac{2 \, \Gamma(j - 1)}{\Gamma(j + k)} = 
\frac{2(2k - 1)\Gamma(k - 1)}{k \, \Gamma(2k)}.$$
The second sum is simplified to
$$\sum_{j = k}^{\infty} \frac{k \, \Gamma(j - 1)}{j \, \Gamma(j + 
	k))} = \frac{\Gamma(k - 1)\Hypergeometric{3}{2}{1, k - 1, k}{2k, 
	k + 
		1}{1}}{\Gamma(2k)}.$$
where $\Hypergeometric{3}{2}{\cdot}{\cdot}{\cdot}$ is a {\em 
	generalized hypergeometric function}. Putting them together, we 
	thus 
have
$$\hat{C}(\infty) \approx \frac{k - 1}{2k - 1}\left(\frac{2(2k - 
	1)}{k} - \Hypergeometric{3}{2}{1, k - 1, k}{2k, k + 
	1}{1}\right).$$
Although hypergeometric functions cannot be written in closed forms 
in general, we derive the analytical results of $\hat{C}(\infty)$ 
for several small values of $k$, and present them in 
Table~\ref{Table:asymcc}. In particular, the estimated clustering 
coefficient for triangulated RANs (i.e., $k = 3)$ based on our 
calculation is $12 \pi^2 - 353/3 \approx 0.7686$, which is more 
accurate than $46/3 - 36 \log(3/2) \approx 0.7366$~\cite[Equation 
(6)]{zhou2005maximal}, according to a simulation 
experiment ($0.7683$ based on the average of $50$ independent 
samples, each of which is run over $10,000$ iterations).

\begin{table}[h]
	\renewcommand{\arraystretch}{1.5}
	\begin{center}
		\begin{tabular}{|c|c|}
			\hline
			Network index ($k$) & $\hat{C}(\infty)$
			\\ \hline
			3 & $12 \pi^2 - \frac{353}{3}$ 
			\\ \hline 4 & $120 \pi^2 - \frac{2367}{2}$
			\\ \hline
			5 & $\frac{2800}{3} \pi^2 - \frac{138161}{15}$ 
			\\ \hline 6 & $6300 \pi^2 - \frac{746131}{12}$
			\\ \hline
			7 & $38808 \pi^2 - \frac{134056533}{350}$ 
			\\ \hline 8 & $224224 \pi^2 - \frac{663900367}{300}$
			\\ \hline
			9 & $1235520 \pi^2 - \frac{26887974331}{2205}$ 
			\\ \hline
			10 & $6563700 \pi^2 - \frac{253941996039}{3920}$
			\\ \hline
		\end{tabular}
	\end{center}
	\caption{Asymptotic clustering coefficients of HDRANs with small 
		indicies $k$}
	\label{Table:asymcc}
\end{table}

\section{Sparsity}
\label{Sec:sparsity}
{\em Sparsity} is a property of common interest in network 
modeling~\cite{singh2015finding, verzelen2015community, 
	vinciotti2013robust}, as well as in 
data 
analytics~\cite{arnold2010specifying, buluc2011implementing}. As 
opposed to 
``dense,'' this 
topology plays a key role when one defines sparse networks. Sparse 
networks have fewer links than the maximum possible number of links 
in the (complete) network of same order. In computer science, sparse 
networks are considered to be somewhere dense or nowhere dense. The 
investigation of sparsity of HDRANs is inspired by an article 
recently published in the American Physics 
Society~\cite{delgenio2011all}. It 
was analytically and numerically proven in the article that the 
probability of a scale-free network being dense is $0$, given that 
the power-law coefficient falls between~$0$ and $2$.

One of the most commonly-used network topology to measure the 
sparsity of a network $G(V, E)$ is the {\em link density} (also 
known as {\em edge density} in the literature): 
$${\rm density}(G) = \frac{|E|}{\binom{|V|}{2}}.$$
For a HDRAN of index $k$, denoted $\apol_{n}^{(k)}$, its link 
density at time $n$ is a decreasing function of $n$, viz.,
$${\rm density}\left(\apol_{n}^{(k)}\right) = 
\frac{E_n^{(k)}}{\binom{V_n^{(k)}}{2}} = \frac{k + nk}{\binom{k + 
		n}{2}} = \frac{2(k + nk)}{(k + n)(k + n - 1)}.$$
Observing that the link density of an HDRAN in any form is 
deterministic given $k$ and $n$, we assert that this topology indeed 
fails to expose the randomness or to capture the structure of 
HDRANs. Other topologies that have been proposed to measure the 
sparsity of both nonrandom and random networks include degeneracy, 
arboricity, maximum average degree, etc. We refer the interested 
readers to~\cite{nesetril2012sparsity} for textbook style 
expositions of these 
topologies and their properties.

In this section, we measure the sparsity of HDRANs via a classical 
metric---the {\em Gini index}~\cite{gini1921measurement}. The Gini 
index which 
appears more often in economics is commonly used to measure the 
inequality of income or wealth~\cite{dalton1920themeasurement, 
	gini1921measurement}. The utilization 
of the Gini index as a sparsity measurement originates in electrical 
engineering~\cite{hurley2009comparing}. More often, the Gini index 
was used to 
evaluate regularity of graphs~\cite{balaji2017thegini, 
	domicolo2020degree, zhang2019thedegree}. The Gini 
index debuted as a sparsity measurement of networks 
in~\cite{goswami2018sparsity}.

\begin{figure}[tbh]
	\centering
	\begin{tikzpicture}[scale=5]
		\draw
		(0, 0) -- (0, 1)
		(0, 0) -- (1, 0)
		(0, 1) -- (1, 1)
		(1, 0) -- (1, 1)
		(0, 0) -- (1, 1);
		\draw[ultra thick]
		(0, 0) to [bend right = 40] (1, 1);
		\draw[fill=gray!20]  
		(0, 0) to [bend right = 40] (1, 1) -- (1, 0) -- (0, 0) -- 
		cycle;
		\node[text width = 0.1cm] at (0.57, 0.43) {$A$};
		\node[text width = 0.1cm] at (0.8, 0.2) {$B$};    
	\end{tikzpicture}
	\label{Fig:exLorenz}
	\caption{An example of typical Lorenz curve}
\end{figure}
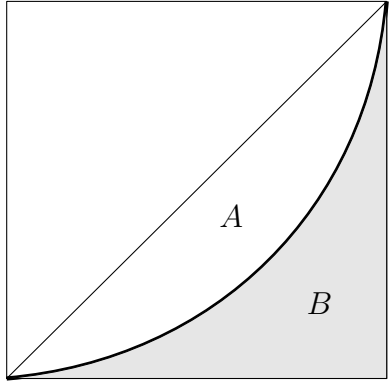

A graphical interpretation of the Gini index is the {\em Lorenz 
	curve}. As portrayed in Figure~\ref{Fig:exLorenz}, the Lorenz 
	curve 
(thick black curve) splits the lower triangle of a unit square into 
$A$ and $B$. 
A well-established relationship between the Gini index and the 
Lorenz curve is that the Gini index of the associated Lorenz curve 
is equal to the ratio of ${\rm Area}(A)$ and ${\rm Area}(A + B)$, 
equivalent to $1 - 2 \times {\rm Area}(B)$.

We construct the Gini index of HDRANs based on vertex degrees. At 
time $n$, there is a total of $k + n$ vertices in $\apol_n^{(k)}$, 
and the {\em admissible} degree set is $\mathcal{J} = \{k, k + 1, 
\ldots, k + n\}$. According to Theorem~\ref{Thm:L1bound}, the mean 
of the proportion of the number of vertices having degree $j \in 
\mathcal{J}$ can be approximated by $\bigl(\Gamma(j) \Gamma(2k - 
1)\bigr)/\bigl(\Gamma(j + k) \Gamma(k - 1)\bigr)$, when $n$ is 
large. For simplicity, let us denote this mean proportion for each 
pair of $j$ and $k$ by $\gamma(j, k)$. These $\gamma(j, k)$'s 
altogether naturally form the Lorenz curve after being rearranged in 
an ascending order. Note that
$$\frac{\partial}{\partial j} \, \gamma(j, k) = \frac{\bigl( 
	(\Psi(j) - \Psi(j + k) \bigr) 2^{2k - 2} (k - 1) \Gamma\left(k - 
	\frac{1}{2}\right) \Gamma(j)}{\Gamma\left(\frac{1}{2}\right) 
	\Gamma(j + k)} < 0,$$ 
where $\Psi(\cdot)$ is the {\em digamma function}, known to be 
increasing on the positive real line. Hence, the function $\gamma(j, 
k)$ is decreasing with respect to $j$.

Specifically, we build the Lorenz curve as follows. The bottom of 
the unit square is equispaced into $(n + 1)$ segments. The bottom 
left vertex is marked $0$ along with vertical value $0$. The 
cumulative proportion value $\sum_{j = k + n - i + 1}^{k + n} 
\gamma(j, k)$ is assigned to the $i$th segmentation from the left. 
There is a total of $n$ segmentations between the bottom left and 
bottom right vertices. Lastly, the vertical value for the bottom 
right vertex is $\sum_{j = k}^{k + n} \gamma(j, k)$. The Lorenz 
curve is comprised by smoothly connecting these assigned values in 
order, from left to right.

In the next lemma, we show that the Lorenz curve that we established 
in the last paragraph is well defined, i.e., the two ends of the 
Lorenz curve respectively coincide with the bottom left and the top 
right corners of the unit square.
\begin{lemma}
	\label{Lem:lorenz}
	We claim that
	$$\lim_{n \to \infty} \sum_{j = k}^{k + n} \frac{\Gamma(j) 
		\Gamma(2k - 1)}{\Gamma(j + k) \Gamma(k - 1)} = 1 \quad 
	\mbox{\textit{and}} \quad
	\lim_{n \to \infty} \frac{\sum_{j = k + n - i + 1}^{k + n} 
		\frac{\Gamma(j) \Gamma(2k - 1)}{\Gamma(j + k) \Gamma(k - 
			1)}}{i/n} = 0.$$
\end{lemma}
The proof of Lemma~\ref{Lem:lorenz} is presented in 
Appendix~\ref{App:lorenz}. Next, we calculate ${\rm Area}(B)$, 
equivalent to integrating the Lorenz curve from $0$ to $1$. For 
large value of $n$, the integration can be approximated by applying 
the {\em trapezoid rule}; that is,
\begingroup
\allowdisplaybreaks
\begin{align*}
	{\rm Area}(B) &\approx \frac{1}{2 (n + 1)} \left[\sum_{j = k + 
		n}^{k + n} \frac{\Gamma(j) \Gamma(2k - 1)}{\Gamma(j + k) 
		\Gamma(k - 1)} \right. 
	\\ &\qquad{}+ \left(\sum_{j = k + n}^{k + n} \frac{\Gamma(j) 
		\Gamma(2k - 1)}{\Gamma(j + k) \Gamma(k - 1)} + \sum_{j = k + 
		n - 
		1}^{k + n} \frac{\Gamma(j) \Gamma(2k - 1)}{\Gamma(j + k) 
		\Gamma(k - 1)}\right)
	\\ &\qquad{} + \cdots + \left.\left(\sum_{j = k + 1}^{k + n} 
	\frac{\Gamma(j) \Gamma(2k - 1)}{\Gamma(j + k) \Gamma(k - 1)} + 
	\sum_{j = k}^{k + n} \frac{\Gamma(j) \Gamma(2k - 1)}{\Gamma(j + 
		k) \Gamma(k - 1)} \right)\right]
	\\&= \frac{1}{2(n + 1)}\left(\frac{3k - 2}{k - 2} - \frac{2^{2k 
			- 1}\bigl((k - 1)n + 2\bigr) \Gamma\left(k - 
		\frac{1}{2}\right)\Gamma(k + n + 1)}{(k - 
		2)\Gamma\left(\frac{1}{2}\right)\Gamma(2k + n)}\right).
	\\&\sim n^{-1} - n^{1 - k}, 
\end{align*}
\endgroup
In what follows, the Gini index of an HDRAN of index $k$ at time $n$ 
is given by
\begin{align*}
	&{\rm Gini}\left(\apol_{n}^{(k)}\right) = 1 - 2 \times {\rm 
		Area}(B)
	\\ &\quad= 1 - \frac{1}{(n + 1)}\left(\frac{3k - 2}{k - 2} - 
	\frac{2^{2k - 1}\bigl((k - 1)n + 2\bigr) \Gamma\left(k - 
		\frac{1}{2}\right)\Gamma(k + n + 1)}{(k - 
		2)\Gamma\left(\frac{1}{2}\right)\Gamma(2k + n)}\right),
\end{align*}
the asymptotic equivalent of which is equal to $1$. A large value of 
Gini index (ranging from $0$ to $1$) indicates an extremely 
nonuniform distribution of vertex degrees, implying that all vertex 
degrees are dominated by only a few classes, whereas a small value 
of Gini index suggests vertex degrees are evenly distributed in 
different degree classes. Thus, we conclude (asymptotically) high 
sparseness of HDRANs.

\begin{figure}[tbh]
	\centering
	\includegraphics[width=\textwidth]{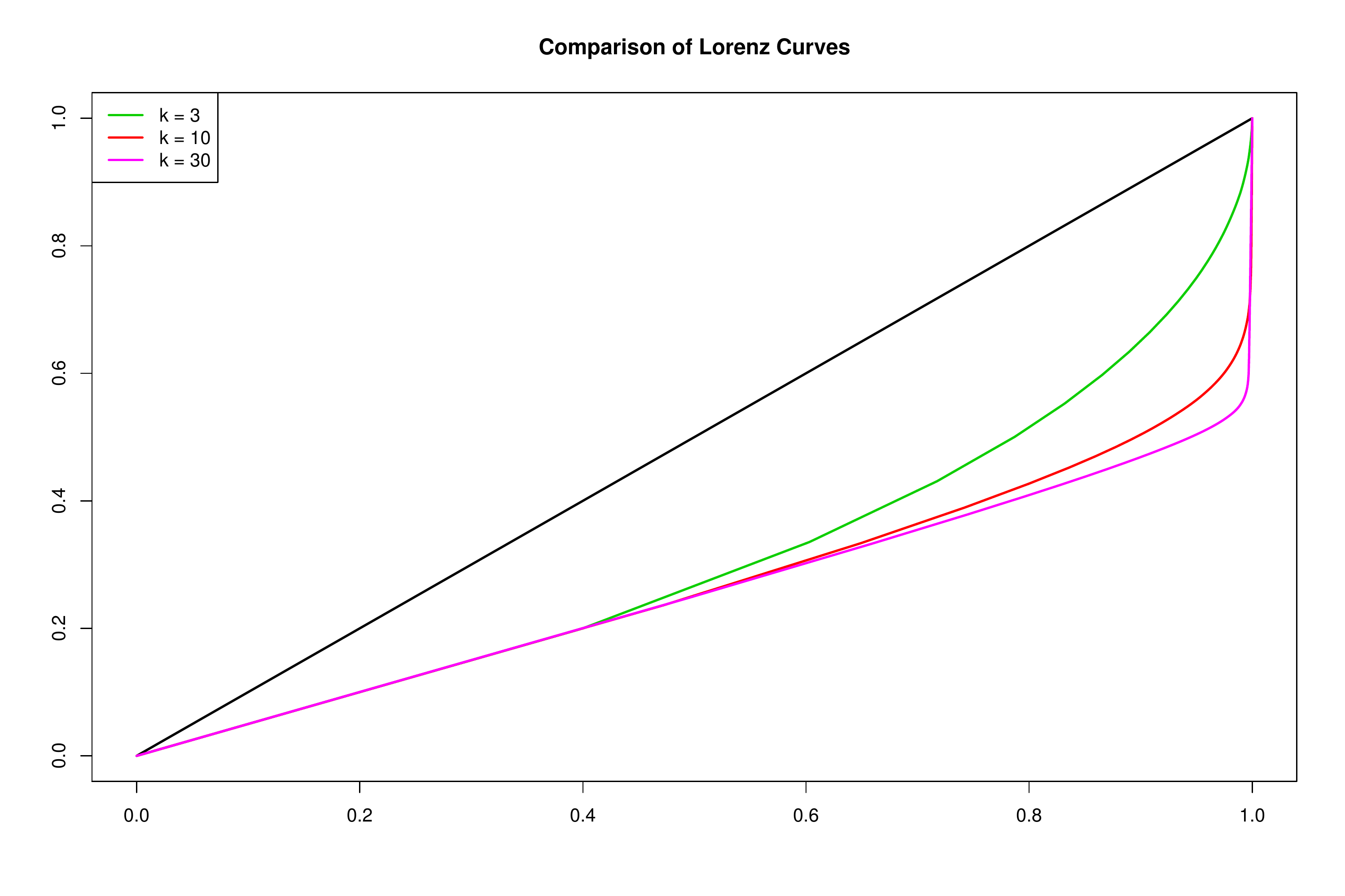}
	\caption{Comparison of Lorenz Curves for simulated HDRANs of $k 
		= 3, 10, 30$ at time $n = 5000$}
	\label{Fig:Lorenz}
\end{figure}

We further verify our conclusion by conducting some simulation 
experiments. In general, each network $G(V, E)$ is associated with a 
unique $|V| \times |V|$ {\em adjacency matrix}, denoted $\matA = 
(A_{ij})$, in which $A_{ij} = 1$ only when there is an edge linking 
vertices $i$ and $j$, for $i, j \in V$; $0$, otherwise. If $G$ is 
undirected, $\matA$ is symmetric. The degree of vertex $i$ thus can 
be represented by the sum of $i$th row or the $i$th column in 
$\matA$, allowing us to compute the Gini index of each simulated 
network~$G$ through $\matA$ accordingly. For each $k = 3, 10, 30$, 
we generate $100$ independent HDRANs at time $n = 5000$. The 
comparison of Lorenz curves (based on the average of cumulative 
degree proportion sequences) is given in Figure~\ref{Fig:Lorenz}.

Besides, we calculate the Gini index of each of the $100$ simulated 
HDRANs (of $k = 3, 10, 30$) at time $50,000$, and take the average; 
The estimated Gini indices are $0.9970330$ (for $k = 3$), 
$0.9990327$ (for $k = 10$), and $0.9997262$ (for $k = 30$). We do 
not show the corresponding Lorenz curves as they are not visually 
distinguishable.

\section{Depth, diameter and distance}
\label{Sec:diameter}
\begin{figure}[tbh]
	\begin{center}
		\begin{tikzpicture}[scale=3.77]
			\draw
			(2.75,0) node [circle,draw] {$0^{(2)}$}
			(3.25,0) node [circle=0.1,draw] {$0^{(1)}$}
			(3,0.866) node [circle=0.1,draw] {$0^{(4)}$}
			(2.55, 0.52) node [circle=0.1,draw] {$0^{(3)}$}
			(3.45, 0.52) node [circle=0.1,draw] {$0^{(5)}$}
			(3.35, 0.866) node [circle,draw] {$\; \; 1 \; \;$}
			(3.55, 0.2) node [circle=0.1,draw] {$\; \; 2 \; \;$};
			\draw
			(2.867,0) -- (3.133,0)
			(2.67,0.52) -- (3.33,0.52)
			(2.842,0.085) -- (2.97, 0.755)
			(3.158, 0.085) -- (3.03, 0.755)
			(2.59, 0.405) -- (2.67, 0.09)
			(3.41, 0.405) -- (3.33, 0.09)
			(2.643, 0.448) -- (3.141, 0.036)
			(3.357, 0.448) -- (2.859, 0.036)
			(2.643, 0.6) -- (2.923, 0.78)
			(3.357, 0.6) -- (3.077, 0.78)
			(3.125, 0.866) to[bend left=15] (3.235, 0.866)
			(2.643, 0.6) to[bend left=75] (3.3, 0.97)
			(2.842, 0.085) to[bend left=45] (3.25, 0.8)
			(3.158, 0.085) to[bend left=15] (3.3, 0.75);
			\draw
			(3.4, 0.63) to[bend right=30] (3.4, 0.75)
			(3.6, 0.32) to [bend right=45] (3.43, 0.78)
			(3.52, 0.32) to [bend left=30] (3.05, 0.765)
			(3.44, 0.25) to [bend left=15] (2.66,0.49)
			(3.436, 0.195) to [bend right=10] (2.842, 0.085)
			(3.47, 0.115) to [bend left=15] (3.37, 0.02);
		\end{tikzpicture}
		\caption{An example of a HDRAN of index 5 at step 2.}
		\label{Fig:RANdepth}
	\end{center} 
\end{figure}
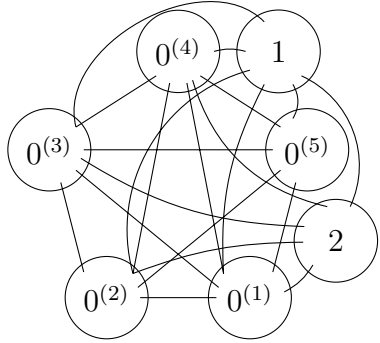

In this section, we investigate several distance-based properties of 
HDRANs. The first measure that we look into is clique-based---{\em 
	depth}, which is defined (for HDRANs) recursively as follows. At 
time $1$, the original $k$-clique is divided into $k$ simplexes, and 
then is deactivated. The depth of each of the active $k$-cliques 
equals $1$. At time $n > 1$, an existing active clique $\clique^{*}$ 
is chosen uniformly at random, and subdivided into $k$ new cliques 
$\clique_1, \clique_2, \ldots, \clique_k$. Then, we have
$${\rm depth}(\clique_i) = {\rm depth}(\clique^{*}) + 1,$$
for all $i = 1, 2, \ldots, k$. An explanatory example of a RAN of 
index $m = 5$ is shown in Figure~\ref{Fig:RANdepth}, where the 
(active) cliques denoted by $\left(1, 0^{(1)}, 0^{(2)}, 0^{(3)}, 
0^{(5)}\right)$, 
$\left(1, 0^{(1)}, 0^{(2)}, 0^{(4)}, 0^{(5)}\right)$, $\left(1, 
0^{(1)}, 0^{(3)}, 0^{(4)}, 0^{(5)}\right)$, and $\left(1, 0^{(2)}, 
0^{(3)}, 0^{(4)}, 0^{(5)}\right)$ have depth $1$; all the rest have 
depth $2$. 

In contrast, {\em distance}, also known as {\em geodesic distance}, 
is a property based on pairwise vertices. In a given network $G(V, 
E)$, the distance between a pair of arbitrary vertices $i, j \in V$, 
denoted $d(i, j)$, is the number of edges in the shortest path (or 
one of the shortest paths) connecting $i$ and $j$. A related 
property, {\em diameter} of network $G$, denoted ${\rm 
	diameter}(G)$, is defined in a max-min manner: the greatest 
	length 
of the shortest paths between every two verticies in $G$, i.e., 
$\max_{i, j \in V} \left\{d(i, j)\right\}$. see~\cite[page 
82]{bondy2008graph} for fundamental properties of the diameter of a 
graph. 
For instance, the diameter of the HDRAN given in 
Figure~\ref{Fig:RANdepth} is $2$, referring to the distance between 
the vertices respectively labeled with $2$ and $0^{(5)}$.

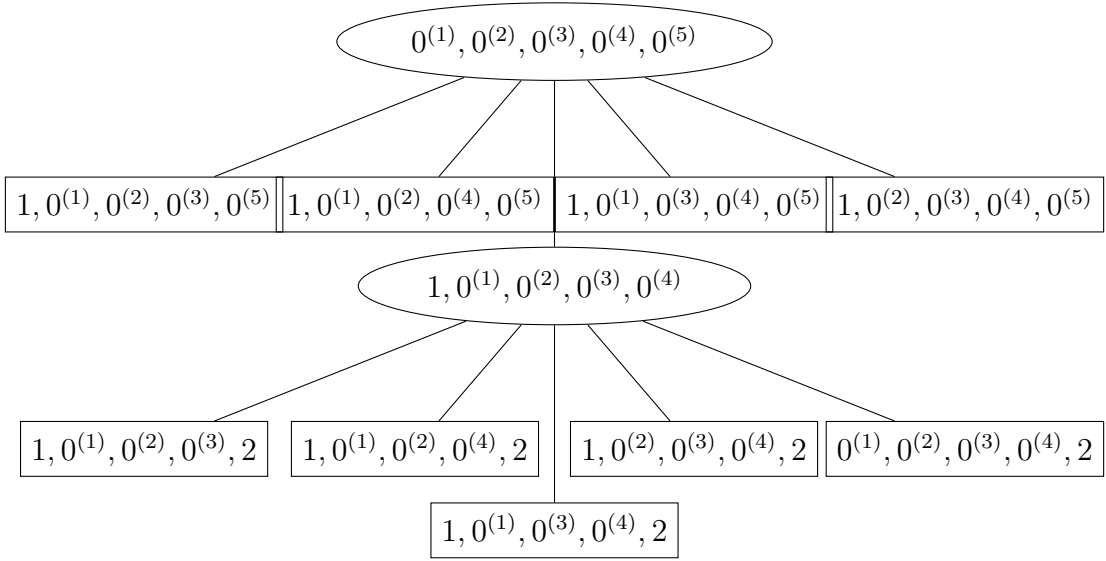
\begin{figure}[tbh]
	\begin{center}
		\centering
		\begin{tikzpicture}[scale = 1.08]
			\node[ellipse, draw] at (0, 0) (v0){$0^{(1)}, 0^{(2)}, 
				0^{(3)}, 0^{(4)}, 0^{(5)}$};
			\node[rectangle, draw] at (-5, -2) (v11){$1, 0^{(1)}, 
				0^{(2)}, 0^{(3)}, 0^{(5)}$};
			\node[rectangle, draw] at (-1.7, -2) (v12){$1, 0^{(1)}, 
				0^{(2)}, 0^{(4)}, 0^{(5)}$};
			\node[ellipse, draw] at (0, -3) (v13){$1, 0^{(1)}, 
				0^{(2)}, 0^{(3)}, 0^{(4)}$};
			\node[rectangle, draw] at (1.7, -2) (v14){$1, 0^{(1)}, 
				0^{(3)}, 0^{(4)}, 0^{(5)}$};
			\node[rectangle, draw] at (5, -2) (v15){$1, 0^{(2)}, 
				0^{(3)}, 0^{(4)}, 0^{(5)}$};
			\node[rectangle, draw] at (-5, -5) (v21){$1, 0^{(1)}, 
				0^{(2)}, 0^{(3)}, 2$};
			\node[rectangle, draw] at (-1.7, -5) (v22){$1, 0^{(1)}, 
				0^{(2)}, 0^{(4)}, 2$};
			\node[rectangle, draw] at (0, -6) (v23){$1, 0^{(1)}, 
				0^{(3)}, 0^{(4)}, 2$};
			\node[rectangle, draw] at (1.7, -5) (v24){$1, 0^{(2)}, 
				0^{(3)}, 0^{(4)}, 2$};
			\node[rectangle, draw] at (5, -5) (v25){$0^{(1)}, 
				0^{(2)}, 0^{(3)}, 0^{(4)}, 2$};
			\draw (v0)--(v11);
			\draw (v0)--(v12);
			\draw (v0)--(v13);
			\draw (v0)--(v14);
			\draw (v0)--(v15);
			\draw (v13)--(v21);
			\draw (v13)--(v22);
			\draw (v13)--(v23);
			\draw (v13)--(v24);
			\draw (v13)--(v25);
		\end{tikzpicture}
		\caption{The evolution of the $5$-ary tree corresponding to 
			that of the HDRAN of index $5$ given in 
			Figure~\ref{Fig:RANdepth}. Elliptic (internal) nodes 
			refer 
			to inactive cliques, whereas rectangular (active) nodes 
			refer to active ones.}
		\label{Fig:karytree}
	\end{center}
\end{figure}

It was introduced in~\cite{darrasse2007degree} that there exists an 
one-to-one 
relation between the evolution of HDRANs (of index $k$) and that of 
$k$-ary trees\footnote{See~\cite[page 224]{storer2002anintroduction} 
	for the 
	definition of $k$-ary tree}. An illustrative example is 
	presented in 
Figure~\ref{Fig:karytree}.  
Active and inactive cliques in HDRANs (of index $k$) respectively 
correspond to external and internal nodes in $k$-ary trees. Thus, 
the total depth of active cliques in $\apol_{n}^{(k)}$ is equivalent 
to the total depth\footnote{In tree structure, the depth of a node 
	is the number of links between the node and the root (of the 
	tree).} 
of external nodes in the corresponding $k$-ary tree at time $n$, 
denoted $\ktree{k}{n}$. In the literature, the total depth of 
external nodes in $\ktree{k}{n}$ is also known as the {total 
	external path}, denoted by $\ext{k}{n}$ in our manuscript. For 
uniformity, we use $\ext{k}{n}$ as the notation for the total depth 
of active cliques in $\apol_{n}^{(k)}$ as well.
\begin{prop}
	\label{Thm:ext}
	Let $\ext{k}{n}$ be the total depth of active cliques in a HDRAN 
	of index $k$ at time $n$. The first two moments of $\ext{k}{n}$ 
	are given by 
	\begin{align*}
		\E\left[\ext{k}{n}\right] &= (kn - n + 1) \sum_{i = 0}^{n - 
			1} \frac{k}{k + (k - 1)i}.
		\\ \E\left[\left(\ext{k}{n}\right)^2\right] &= \bigl((k - 
		1)n + m\bigr) \bigl((k - 1)n + 1\bigr) k E(k, n) + 
		O\left(n^2 \log{n}\right),
	\end{align*}
	where $E(k, n)$ is a function of $k$ and $n$, given in 
	Appendix~\ref{App:ext}.
\end{prop}
The proof of Proposition~\ref{Thm:ext} also can be found in 
Appendix~\ref{App:ext}. As we know that 
$$\sum_{i = 0}^{n - 1} \frac{k}{k + (k - 1)i} \sim \frac{k}{k - 1} 
\log{n},$$
for large $n$, we hence conclude that the leading order of the 
asymptotic expectation of $\ext{k}{n}$ is $kn\log{n}$.

The diameter of HDRANs is also considered. In~\cite{frieze2014some}, 
the 
authors established an upper bound for the diameter of planar RANs 
by utilizing a known result of the height of weighted $k$-ary 
trees~\cite[Theorem 5]{broutin2006large}, i.e.,
$${\rm diameter} \left(\apol_{n}^{(3)}\right) \le \rho \log n,$$
where $\rho = 1/\eta$, and $\eta$ is the unique solution greater 
than $1$ for $\eta - 1 - \log \eta = \log 3$. This upper bound can 
be extended to $\apol_{n}^{(k)}$ effortlessly; that is, 
$${\rm diameter} \left(\apol_{n}^{(k)}\right) \le \frac{2}{\rho^{*} 
	(k - 1)} \log n,$$
where $\rho^{*} = 1 / \eta^{*}$ is the unique solution greater than 
$1$ for $\eta^{*} - 1 - \log \eta^{*} = \log k$. In addition, the 
authors of~\cite{ebrahimzadeh2014onlongest} proved ${\rm diameter} 
\left(\apol_{n}^{(3)}\right) \overset{a.s.}{\sim} c \log n$ by 
estimating the height of a class of specifically-designed random 
trees. The value of $c$ is approximately $1.668$. The asymptotic 
expression of the diameter of more general $\apol_{n}^{(k)}$ was 
developed by~\cite{cooper2014theheight} and 
by~\cite{kolossvary2016degrees}. The 
approach in~\cite{cooper2014theheight} was to utilize known results 
of 
continuous-time branching processes coupled with recurrence methods, 
and the authors of~\cite{kolossvary2016degrees} coped with 
difficulties by 
characterizing vertex generations. We only state (without repeating 
the proof) the weak law of the diameter of $\apol_{n}^{(k)}$ 
from~\cite{cooper2014theheight} (with a minor tweak) in the next 
theorem.
\begin{theorem}[\mbox{~\cite[Theorem 2]{cooper2014theheight}}]
	For $k \ge 3$, with high probability, we have
	$${\rm diameter} \left(\apol_{n}^{(k)}\right) \sim c \log n,$$
	where $c$ is the solution of 
	$$\frac{1}{c} = \sum_{\ell = 0}^{k - 1} \frac{k - 1}{\ell + a(k 
		- 1)},$$
	in which the value of $a$ is given by
	$$\frac{\Gamma(k + 1) \Gamma(ka)}{\Gamma\bigl((k - 1)a + 
		k\bigr)} \exp\left\{ \sum_{\ell = 0}^{k - 1} \frac{(k - 1)(a 
		+ 
		1) - 1}{\ell + (k - 1)a} \right\} = 1.$$
	Especially, as $k \to \infty$,
	$$c \sim \frac{1}{k \log 2}.$$
\end{theorem}

A topological measure related to distance is the {\em Wiener index}, 
which was proposed by the chemist Harry 
Wiener~\cite{wiener1947structural} to 
study molecular branching of chemical compounds. For a network $G(V, 
E)$, the Wiener index is defined as the sum of distances of all 
paired vertices, i.e., $W(G) = \sum_{i, j \in V} d(i, j)$. The 
Wiener index has been extensively studied for random 
trees~\cite{dobrynin2001wiener, neininger2002thewiener}. For other 
random structures, we 
refer the readers to~\cite{bereg2007wiener, fuchs2015thewiener, 
	janson2003thewiener, wagner2006aclass, wagner2007ontheaverage, 
	wagner2012onthewiener}. 

The methodologies for computing the Wiener index of random trees, 
however, are 
not adaptable to the study of RANs, as the bijection between RANs 
and $k$-ary trees is based on a clique-to-node mapping. The high 
dependency of active cliques (sharing vertices and edges) 
substantially increases the challenge of formulating mathematical 
relation between distance (vertex-based) and depth (clique-based).

There is only a few articles studying distance or related properties 
in RANs. In~\cite{kolossvary2016degrees}, the authors proved that 
the distance 
of two arbitrary vertices in a HDRAN has both mean and variance of 
order $\log n$, and that this distance follows a Gaussian law 
asymptotically. However, it seems difficult to extend this result to 
the Wiener index, as the covariance structure of the distances (of 
all paired vertices) is unspecified. Planar RANs ($\apol_{n}^{(3)}$) 
were considered in~\cite{bodini2008distances}. In this article, the 
dominant term 
of the total distance of all pairs of vertices was shown to be 
$\sqrt{3 \pi} n^{5/2}/22$. The main idea was to consider an 
enumerative generating function of the total distance, and then 
decompose the total distance into interdistance and extradistance. 
This approach can be extended to HDRANs of small network index~$k$, 
but seemingly not applicable to HDRANs with general index $k$. 
Therefore, the Wiener index of HDRANs remains an open problem.

\begin{figure}[tbh]
	\centering
	\begin{minipage}{0.48\textwidth}
		\includegraphics[scale=0.23]{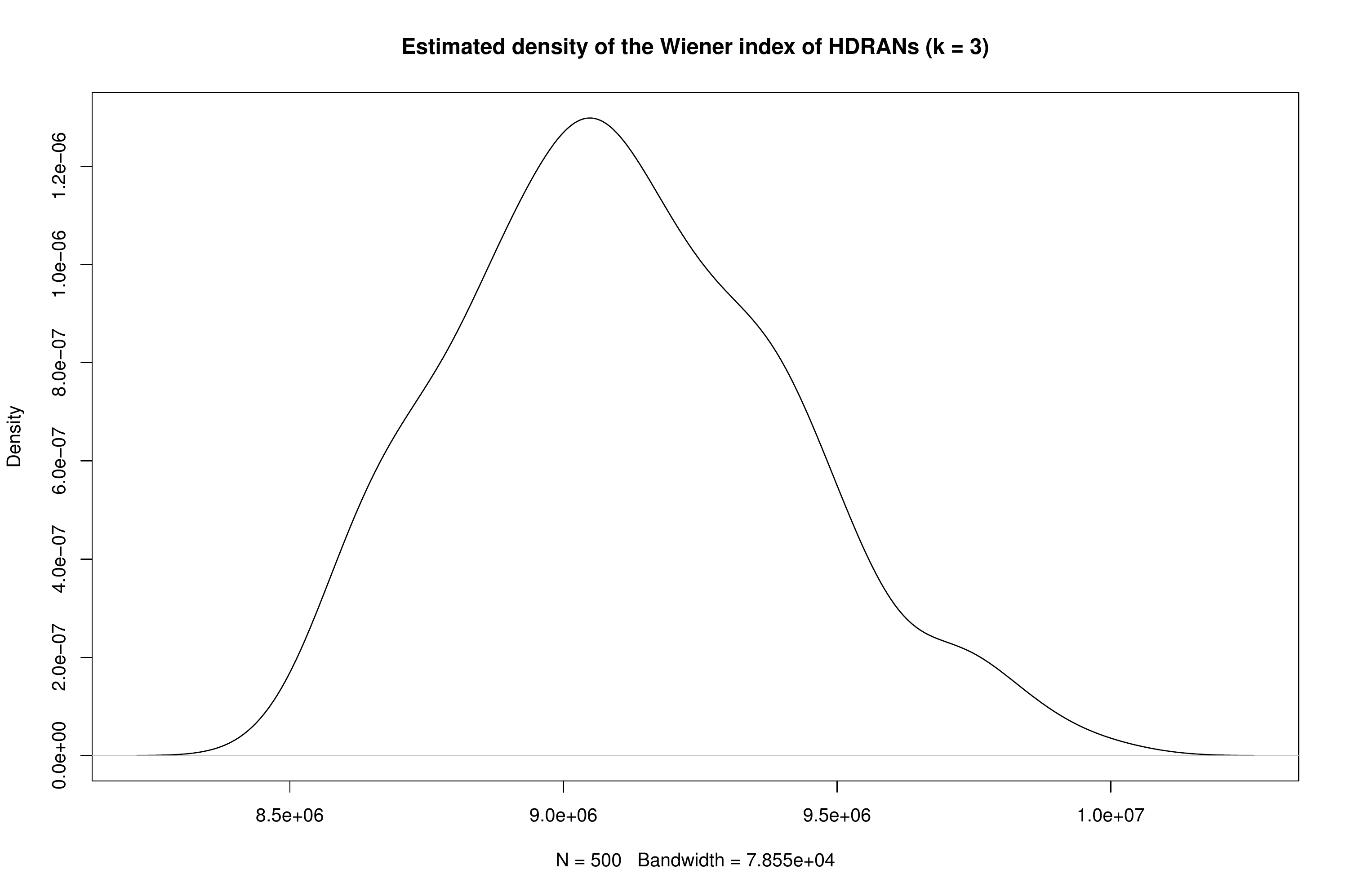}
	\end{minipage}
	\begin{minipage}{0.48\textwidth}
		\includegraphics[scale=0.23]{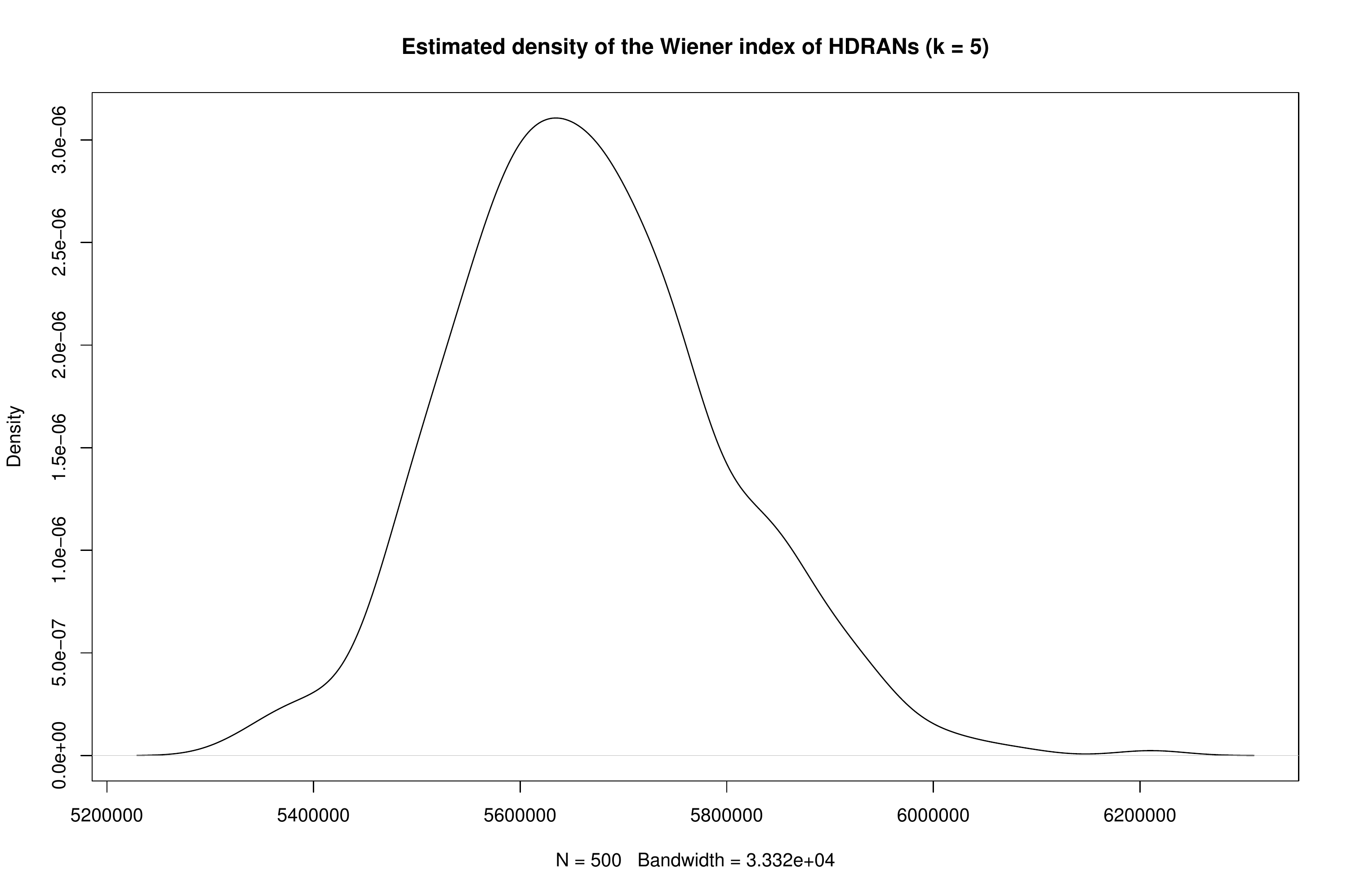}
	\end{minipage}
	\\
	\begin{minipage}{0.48\textwidth}
		\includegraphics[scale=0.23]{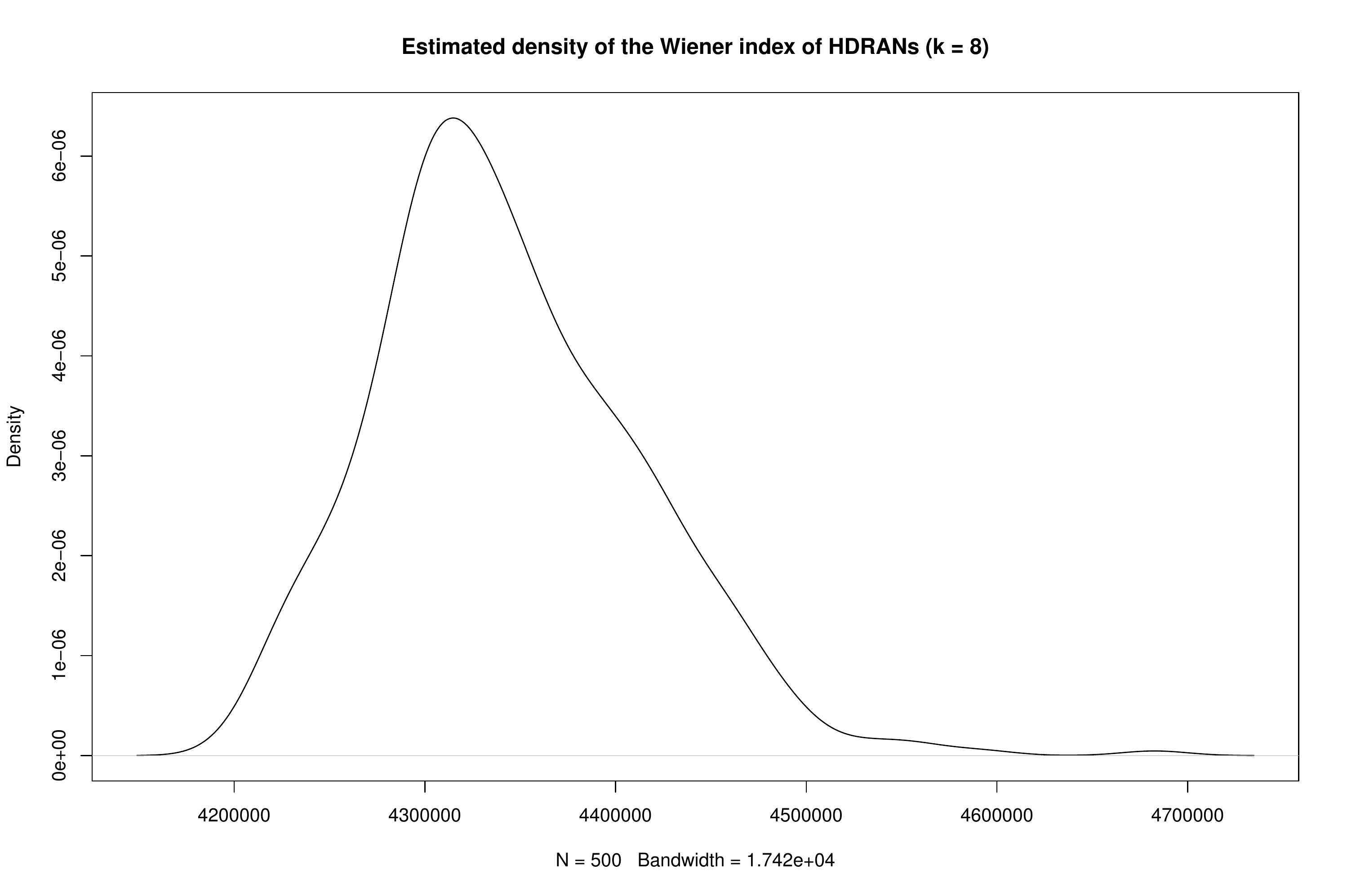}
	\end{minipage}
	\begin{minipage}{0.48\textwidth}
		\includegraphics[scale=0.23]{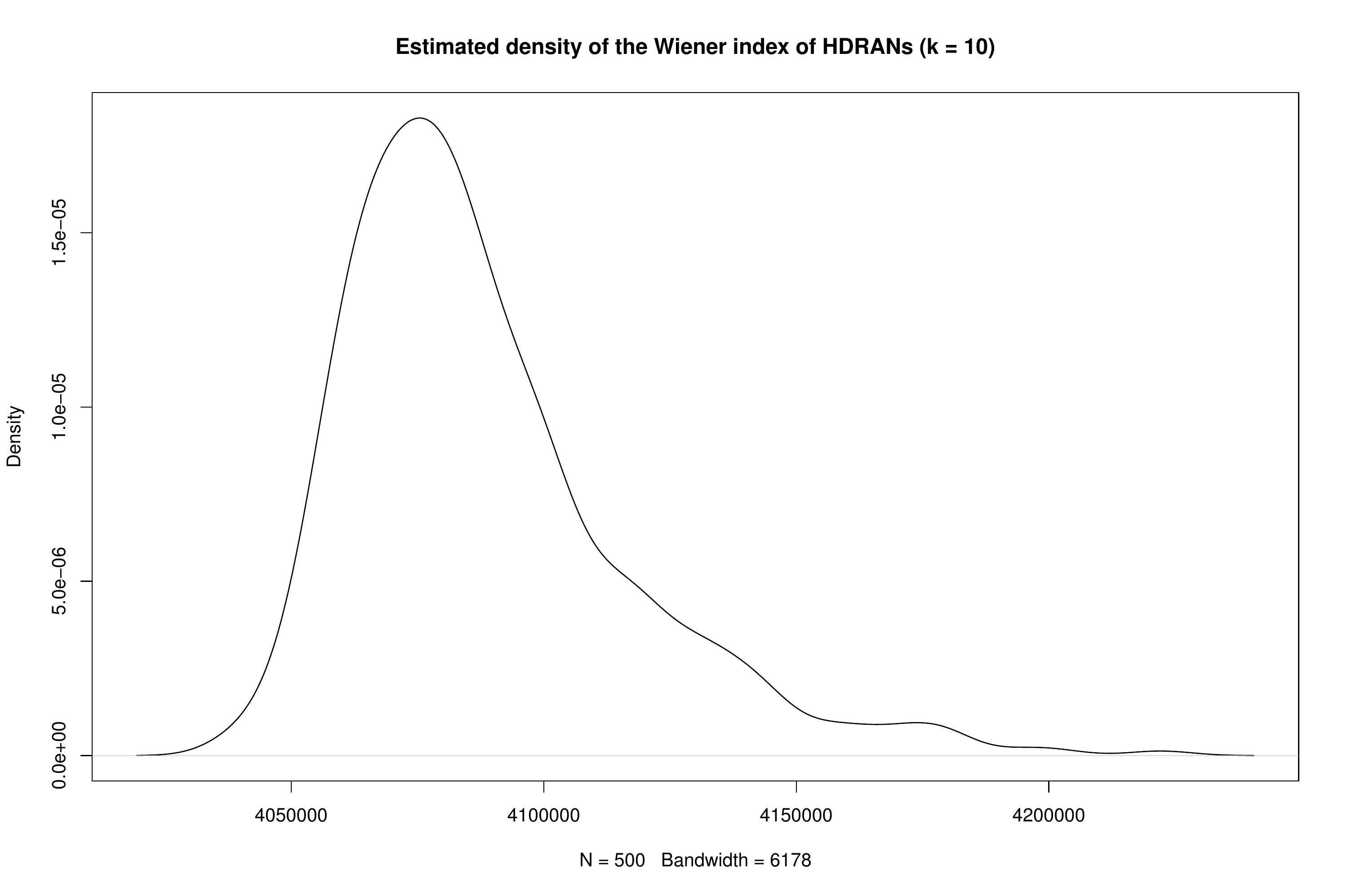}
	\end{minipage}
	\caption{Density estimation of the Wiener indices of HDRANs for 
		$k = 3, 5, 8, 10$}
	\label{Fig:densityest}
\end{figure}

We numerically explore the Wiener index of HDRANs via a series of 
simulations. For $k = 3, 5, 8, 10$, we generate $500$ independent 
HDRANs at time $2,000$, calculate the Wiener index for each 
simulated 
HDRAN, and use the kernel method to estimate the density. The plots 
of the estimated densities are presented in 
Figure~\ref{Fig:densityest}, where we find that they are 
approximately bell-shaped, but not symmetric (positively skewed). By 
observing these patterns, we conjecture that the limiting 
distribution of the Wiener index of HDRANs does not follow a 
Gaussian law. 
In addition, for each $k$, we apply the {\em Shapiro-Wilk} test to 
the simulated data comprising $500$ Wiener indices, and receive the 
following $p$-values: $0.0003$ for $k = 3$; $0.0024$ for $k = 5$; 
$9.56 \times 10^{-8}$ for $k = 8$; and $0$ for $k = 10$. These 
$p$-values are all statistically significant, in support of our 
conjecture.

\section{Concluding remarks}
\label{Sec:concluding}
Finally, we address some concluding remarks and propose some 
future work. We investigate several properties of high-dimensional 
random Apollonian networks in this paper. Two types of degree 
profiles are considered. For the first type, we show that the number 
of vertices of a given degree concentrates on its expectation with 
high probability. In the proof of Theorem~\ref{Thm:degreedist}, we 
derive the $L_1$ limit of $X_{n, j}^{(k)}$, i.e.,
$$\lim_{n \to \infty} \E \left[X_{n, j}^{(k)} \right]= 
\frac{\Gamma(j) \Gamma(2k - 1)}{\Gamma(j + k) \Gamma(k - 1)} \, n,$$
which suggests that the asymptotic expectation of $X_{n, j}^{(k)}$ 
experiences a phase transition. There are two regimes. According to 
the Stirling's Approximation, we have
$$\E \left[X_{n, j}^{(k)} \right] \sim
\begin{cases}
	\frac{\Gamma(j) \Gamma(2k - 1)}{\Gamma(j + k) \Gamma(k - 1)} \, 
	n, &\qquad \mbox{for fixed }j;
	\\ \frac{\Gamma(2k - 1)}{\Gamma(k - 1)} \frac{n}{j^k}, &\qquad 
	\mbox{for }j \to \infty,
\end{cases}
$$
as $n \to \infty$.

For the second type of degree profile, the degree of a vertex of a 
given label, we develop the probability mass function and the exact 
moments by applying the analytic combinatorics methods and the 
results in triangular \polya\ urns.

The next two properties that we investigate are the small world 
property measured by the local clustering coefficient and the 
sparsity measured by a proposed Gini index. We conclude that HDRANs 
are highly clustered and sparse.

The last several properties that we look into are based on node 
distance. According to an one-to-one relation between HDRANs and 
$k$-ary trees, we compute the first two moments of the total depth 
of active cliques in HDRANs. We also numerically study the Wiener 
index, and conjecture that its limiting distribution is not normal 
based on simulation results. The diameter of HDRANs is retrieved 
from~\cite{cooper2015long}.

To end this section, we propose some future work. Our conjecture of 
non-normality of the Wiener index is based on numerical experiments. 
A more rigorous proof is needed. There remain many open problems for 
HDRANs, such as the length of the longest path and the highest 
vertex degree. One suggests carrying out some studies of the 
stochastic 
processes that 
take place on HDRANs, especially those with applications to 
mathematical physics, such as percolation and diffusion. We will 
look into these open problems and report the results elsewhere.

\newpage

\appendix
\section{Proof of Theorem~\ref{Thm:L1bound}}
\label{App:L1bound}
To establish a sharp bound for the difference between the 
expectation of $X_{n, j}^{(k)}$ and its $L_1$ limit (after being 
properly scaled) for $j \ge k$, we distinguish the case of $j = k$ 
and the case of $j > k$.
\subsection{The case of $j = k$}
The vertices of degree $k$ form a special class in RANs of index 
$k$---{\em terminal vertices}. Terminal vertices never recruit 
newcomers since their first appearance in the network. A stronger 
almost sure limit of $X_{n, k}^{(k)}$ was developed via a two-color 
{\em \polya\ urn model}\footnote{We refer the interested readers 
	to~\cite{mahmoud2009polya} for text style exposition of \polya\ 
	urn 
	models.} 
in~\cite{zhang2016thedegree}; that is,
$$\frac{X_{n, k}^{(k)}}{n} \convas \frac{k - 1}{2k - 1}.$$
According to the result in Theorem~\ref{Thm:L1bound}, the $L_1$ 
bound for $X_{n, k}^{(k)}$ is given by
\begin{equation}
	\label{Eq:degreeL1k}
	\left|\E \left[X_{n, k}^{(k)}\right] - \frac{k - 1}{2k - 1} 
	n\right| \le \frac{2k^2}{2k - 1}.
\end{equation}
We prove this result via an induction on $n \in \mathbb{N}$. 
Obviously, Equation~(\ref{Eq:degreeL1k}) is valid for $n = 1$, as we 
have $\E \left[X_{1, k}^{(k)}\right] = X_{1, k}^{(k)} = 1$. Assume 
that Equation~(\ref{Eq:degreeL1k}) holds up to some integer $m$. 
Consider a recursive relation between $X_{n + 1, k}^{(k)}$ and 
$X_{n, k}^{(k)}$ for all $n \ge 1$. Denote ${\rm deg}_v(n)$ the 
degree of vertex $v$ at time $n$, and let $\indicator \{{\rm 
	deg}_v(n) = j\}$ be an indicator function which equals $1$, if 
	${\rm 
	deg}_v(n) = j$; $0$, otherwise. In general, the expected value 
	of 
$X_{n, j}^{(k)}$ can be written in terms of the following:
\begin{equation}
	\label{Eq:indicatorrelation}
	\E \left[X_{n, j}^{(k)}\right] = \sum_{v} \E\left[\indicator 
	\{{\rm deg}_v(n) = j\}\right].
\end{equation}
Besides, we have the following almost-sure relation of the degree of 
$v$ between time $n + 1$ and $n$:
$$\indicator \{{\rm deg}_v(n + 1) = k\} = \indicator \{{\rm 
	deg}_v(n) = k\} \Prob(v\mbox{ is not chosen at time }n + 1).$$
In what follows, we obtain a recurrence between the first moments of 
$X_{n + 1, k}^{(k)}$ and $X_{n, k}^{(k)}$; that is,
$$\E \left[X_{n + 1, k}^{(k)}\right] = \E \left[X_{n, 
	k}^{(k)}\right] \left(1 - \frac{k}{(k - 1)n + 1}\right) + 1.$$
For the inductive step $n = m + 1$, we have
\begin{align*}
	\allowdisplaybreaks
	&\left|\E \left[X_{m + 1, k}^{(k)}\right] - \frac{k - 1}{2k - 1} 
	(m + 1) \right| 
	\\ &\quad{}= \left|\E \left[X_{m, k}^{(k)}\right] \left(1 - 
	\frac{k}{(k - 1)m + 1}\right) + 1 - \frac{k - 1}{2k - 1} (m + 1) 
	\right|
	\\ &\quad{}\le \left|\left( \E \left[X_{m, k}^{(k)}\right] - 
	\frac{k - 1}{2k - 1} m\right) \left(1 - \frac{k}{(k - 1)m + 
		1}\right)\right| 
	\\ &\quad\qquad{}+ \left|\frac{k}{2k - 1} - \frac{k(k - 1) 
		m}{(2k - 1)((k - 1)m + 1)}\right|
	\\ &\quad\le \frac{2k^2}{2k - 1}\left(1 - \frac{k}{(k - 1)m + 
		1}\right) + \frac{k}{(2k - 1)((k - 1)m + 1)},
\end{align*}
which completes the proof.
\subsection{The case of $j > k$}
Due to high dependency of HDRAN structure in the network growth, it 
is difficult to use classical probabilistic methods, such as \polya\ 
urns or recurrence methods, to determine the $L_1$ limit or 
establish an $L_1$ bound for $X_{n, j}^{(k)}$ for general $j > k$. 
The reason is that the recurrence for the expectation of $X_{n, 
	j}^{(k)}$ does not have an analytic solution.

For this case, we prove the theorem by a two-dimensional induction 
on $n = \{1, 2, 3, \ldots\}$ and $j = \{k, k + 1, \ldots, k + n - 
1\}$. Consider an infinite lower triangle table in which the rows 
are indexed by $n$ and the columns are indexed by $j$. A 
illustrative diagram of the inductive progression can be found 
in~\cite[page 69]{zhang2016onproperties}. The leftmost column and 
the diagonal 
of the triangle jointly form the bases of the induction. Notice that 
the leftmost column refers to the case of $j = k$, which has already 
been verified for all $n$. The basis on the diagonal can be proved 
in an analogous manner. We omit the details here. 

Assume that the result stated in the theorem holds up to $j = \ell > 
k$. Before proving the inductive step, we establish a 
two-dimensional recursive relation for $\E \left[X_{n, j}^{(k)} 
\right]$ for $n$ and $j$. Since the degree for each vertex in the 
network increases at most by one at each evolutionary step, we 
observe an almost-sure relation for ${\rm deg}_v(n)$ as follows:
\begingroup
\allowdisplaybreaks
\begin{align*}
	\indicator \{{\rm deg}_v(n + 1) = j\} &= \indicator \{{\rm 
		deg}_v(n) = j\} \Prob(v\mbox{ is not chosen at time }n + 1)
	\\ &\qquad{} + \indicator \{{\rm deg}_v(n) = j - 1\} 
	\Prob(v\mbox{ is chosen at time }n + 1)
	\\ &= \indicator \{{\rm deg}_v(n) = j\} \left(1 - \frac{j}{(k - 
		1)n + 1}\right) 
	\\ &\qquad{}+ \indicator \{{\rm deg}_v(n) = j - 1\} \frac{j - 
		1}{(k - 1)n + 1}.
\end{align*}
\endgroup

According to Equation~(\ref{Eq:indicatorrelation}), we then obtain a 
recurrence for $\E \left[X_{n, j}^{(k)}\right]$; namely,
$$\E \left[X_{n + 1, j}^{(k)}\right] = \E \left[X_{n, 
	j}^{(k)}\right] \left(1 - \frac{j}{(k - 1)n + 1}\right) + \E 
\left[X_{n, j - 1}^{(k)}\right] \frac{j - 1}{(k - 1)n + 1}.$$
To the best of our knowledge, the recurrence above does not have an 
analytic solution. We exploit a well-known result 
in~\cite{chung2006complex},  to 
compute the asymptotic expectation of $\E \left[X_{n, 
	j}^{(k)}\right]$, and to determine the value of $b_{j, k}$ 
subsequently. It was shown in~\cite{chung2006complex} that a 
sequence 
$\{\alpha_n\}$ which satisfies the recurrence
$$\alpha_{n + 1} = \left(1 - \frac{\beta_n}{n + \xi}\right) \alpha_n 
+ \gamma_n$$
for $n \ge n_0$ such that $\lim_{n \to \infty} \beta_n = \beta > 0$ 
and $\lim_{n \to \infty} \gamma_n = \gamma$ has the following 
limiting result:
$$\lim_{n \to \infty} \frac{\alpha_n}{n} = \frac{\gamma}{1 + 
	\beta}.$$
Consider the following settings: $\alpha_n = \E \left[X_{n, 
	j}^{(k)}\right]$, $\beta_n = j/(k - 1)$, $\xi = 1/(k - 1)$, and 
$\gamma_n = \E \left[X_{n, j - 1}^{(k)}\right] \frac{j - 1}{(k - 1)n 
	+ 1}$. We then have
$$\lim_{n \to \infty} \frac{\E \left[X_{n, j}^{(k)}\right]}{n} = 
b_{j - 1, k} \frac{j - 1}{j + k - 1},$$
which in fact establishes a heirachical recurrence for $b_{j, k}$ 
for $j \ge k$; that is,
$$b_{j, k} = b_{j - 1, k} \frac{j - 1}{j + k - 1}.$$
with the initial value $b_{k, k} = (k - 1)/(2k - 1)$. We solve the 
recurrence to get
$$b_{j, k} = \frac{(k - 1)\Gamma(j)\Gamma(2k)}{(2k - 1) \Gamma(j + 
	k) \Gamma(k)} = \frac{\Gamma(j)\Gamma(2k - 1)}{\Gamma(j + k) 
	\Gamma(k - 1)}.$$
We are now at the position to prove the inductive step. For $j = 
\ell + 1$, we have
\begingroup
\allowdisplaybreaks
\begin{align*}
	&\E \left[X_{n + 1, \ell + 1}^{(k)}\right] - b_{\ell + 1, k} (n 
	+ 1)  
	\\&=  \E \left[X_{n, \ell}^{(k)}\right] \left(1 - \frac{\ell}{(k 
		- 1)n + 1}\right) + \E \left[X_{n, \ell - 1}^{(k)}\right] 
	\frac{\ell - 1}{(k - 1)n + 1} - b_{\ell + 1, k} (n + 1)
	\\&= \left(\E \left[X_{n, l}^{(k)}\right] - b_{l, k} n \right) 
	\left(1 - \frac{l}{(k - 1)n + 1}\right) + \E \left[X_{n, l - 
		1}^{(k)}\right] \frac{l - 1}{(k - 1)n + 1} 
	\\&\qquad{}+ b_{l - 1, k} \left(\frac{l - 1}{l + k - 
		1}\right)\left[\left(1 - \frac{l}{(k - 1)n + 1}\right)n - 
	\frac{l}{l + k}(n + 1)\right]
	\\&\le \left(\E \left[X_{n, l}^{(k)}\right] - b_{l, k} n \right) 
	\left(1 - \frac{l}{(k - 1)n + 1}\right) 
	\\&\qquad{}+ \left(\E \left[X_{n, l - 1}^{(k)}\right] - b_{l - 
		1, k} n \right) \frac{l - 1}{(k - 1)n + 1} .
\end{align*}
\endgroup

Therefore, we arrive at 
$$\left|\E \left[X_{n + 1, l + 1}^{(k)}\right] - b_{l, k} (n + 
1)\right| \le \frac{2k^2}{k - 1}\left(1 - \frac{1}{(k - 1)n + 
	1}\right) \le \frac{2k^2}{k - 1},$$
which completes the proof.

\section{Proof of Theorem~\ref{Thm:pbound}}
\label{App:pbound}
At first, we state the {\em Azuma-Hoeffding inequality}: Let 
$(M_n)_{n = 0}^{N}$ be a martingale sequence such that $|M_{n + 1} - 
M_n| \le c$ for all $0 \le n \le N - 1$. For an arbitrary $\lambda > 
0$, we have 
$$\Prob\left[|M_n - M_0| \ge \lambda \right] \le 
e^{-\frac{\lambda^2}{2c^2N}}.$$
Let $(\Omega, \mathbb{F}, \Prob)$ be the probability space induced 
by a HDRAN after $n$ insertions. Fix $j \ge k$, let $(M_{i})_{i = 
	0}^{n}$ be the martingale sequence defined as $M_i = \E 
\left[X^{(k)}_{n, j} \, \Big| \, \F_i\right]$, where $\F_i$ is the 
$\sigma$-field generated by the HDRAN in the first $i$ steps. Recall 
that $M_0 = \E \left[X^{(k)}_{n, j} \, \Big| \, \F_0\right] = \E 
\left[X^{(k)}_{n, j}\right]$ and $M_n = \E \left[X^{(k)}_{n, j} \, 
\Big| \, \F_n\right] = X^{(k)}_{n, j}$. Consider two stochastic 
sequences of choosing cliques, $\mathcal{S} := \clique_1, \clique_2, 
\ldots, \clique_{s - 1}, \clique_s, \ldots$ and 
$\mathcal{S}^{\prime} := \clique_1, \clique_2, \ldots, \clique_{s - 
	1}, \clique^{\prime}_s, \ldots$
in which the first different choices of cliques appear at time $s$. 
The change of the number of vertices of degree $j$ at time $s$ is at 
most $2k$, referring to the $2k$ vertices involved in the cliques 
$\clique_s$ and $\clique^{\prime}_s$. We manipulate the choices of 
cliques in the two sequences after time $s$ as follows:
\begin{itemize}
	\item If an active clique inside of $\clique_s$ in $\mathcal{S}$ 
	is chosen, we select the same clique in $\mathcal{S}^{\prime}$;
	\item If an active clique outside of $\clique_s$ in 
	$\mathcal{S}$ is chosen, we select the a clique inside 
	$\clique_s^{\prime}$ according to a preserved isomorphic mapping 
	in $\mathcal{S}^{\prime}$.
\end{itemize}
Noticing that $\clique_s$ and $\clique^{\prime}_s$ are arbitrary, we 
conclude that the difference between the number of vertices with 
degree $j$ in $\mathcal{S}$ and $\mathcal{S}^{\prime}$ is at most 
$2k$ in average. Thus, we have $|M_{i + 1} - M_i| \le 2k$, for all 
$0 \le i \le n - 1$. The result in Theorem~\ref{Thm:pbound} is then 
immediately obtained by applying the Azuma-Hoeffding inequality.

\section{Proof of Lemma~\ref{Lem:lorenz}}
\label{App:lorenz}
\subsection{The first part}
The identity in the left-hand side of the equation is summable:
\begin{align*}
	\sum_{j = k}^{k + n} \frac{\Gamma(j) \Gamma(2k - 1)}{\Gamma(j + 
		k) \Gamma(k - 1)} &= \frac{\Gamma(2k - 1)}{\Gamma(k - 1)} 
	\sum_{j = k}^{k + n}\frac{\Gamma(j)}{\Gamma(j + k)}
	\\ &= \frac{\Gamma(2k - 1)}{\Gamma(k - 1)} \left(-\frac{\Gamma(k 
		+ n + 1)}{(k - 1)\Gamma(2k + n)} + 
	\frac{\Gamma\left(\frac{1}{2}\right)2^{-2k + 2}(2k - 1)}{2(k - 
		1)\Gamma\left(k + \frac{1}{2}\right)}\right)
	\\ &= -\frac{\Gamma(2k - 1) \Gamma(k + n + 1)}{\Gamma(k) 
		\Gamma(2k + n)} + 1.
\end{align*}
As $n \to \infty$, we apply the Stirling's approximation to get
$$\lim_{n \to \infty} \left(-\frac{\Gamma(2k - 1) \Gamma(k + n + 
	1)}{\Gamma(k) \Gamma(2k + n)} + 1\right) \sim \lim_{n \to 
	\infty} 
\left(n^{1 - k}\right) + 1 = 1,$$
for $k \ge 3$.
\subsection{The second part}
We change the lower bound of the summation to $k + n - i + 1$, and 
get
\begin{align*}
	\sum_{j = k + n - i + 1}^{k + n} \frac{\Gamma(j) \Gamma(2k - 
		1)}{\Gamma(j + k) \Gamma(k - 1)} &= \sum_{j = k}^{k + n} 
	\frac{\Gamma(j) \Gamma(2k - 1)}{\Gamma(j + k) \Gamma(k - 1)} - 
	\sum_{j = k}^{k + n - i} \frac{\Gamma(j) \Gamma(2k - 
		1)}{\Gamma(j + k) \Gamma(k - 1)}
	\\ &= \frac{\Gamma(2k - 1) \Gamma(k + n - i + 1)}{\Gamma(k) 
		\Gamma(2k + n - i)} -\frac{\Gamma(2k - 1) \Gamma(k + n + 
		1)}{\Gamma(k) \Gamma(2k + n)}
	\\ &\sim n^{1 - k},
\end{align*}
by the Stirling's approximation. For $k \ge 3$, we have $n^{1 - k} = 
o\left(n^{-1}\right)$, which completes the proof.

\section{Proof of Proposition~\ref{Thm:ext}}
\label{App:ext}
The proof is based on $\ktree{k}{n}$. A related property of 
$\ext{k}{n}$ is the number of external nodes in $\ktree{k}{n}$, 
denoted by $L_{n}^{(k)}$. According to the evolution of $k$-ary 
trees, we have
$$L_{n}^{(k)} = 1 + n(k - 1).$$
At time $n$, we enumerate all external nodes in $\ktree{k}{n}$ with 
respect to a preserved isomorphic mapping, e.g., from top to bottom, 
and from left to right for the external nodes at the same level. Let 
$D^{(k)}_{i, n}$ be the depth of the external node labeled with $i$ 
in $\ktree{k}{n}$. We thus have 
$$\ext{k}{n} = \sum_{i = 1}^{L_n^{(k)}} D^{(k)}_{i, n}.$$

Note that the quantity of $\ext{m}{n}$ increases monotonically with 
respect to $n$, and the amount of increase depends on the depth of 
the node sampled at each time point. Suppose that the external node 
labeled with $i$ is selected upon time $n$. The increment from 
$\ext{k}{n - 1}$ to $\ext{k}{n}$ is
$$ k\left(D^{(k)}_{i, n - 1} + 1\right) - D^{(k)}_{i, n - 1} = (k - 
1) D^{(k)}_{i, n - 1} + k,$$
leading to the following almost-sure relation conditioning on 
$\ktree{m}{n - 1}$ and label $i$, i.e.,
\begin{equation}
	\label{Eq:extinc}
	\ext{k}{n} = \ext{k}{n - 1} + (k - 1) D^{(k)}_{i, n - 1} + k
\end{equation}
We obtain a recurrence for $\E\left[\ext{k}{n}\right]$ by averaging 
$i$ out in Equation~(\ref{Eq:extinc}); that is,
\begin{align*}
	\E\left[\ext{k}{n} \Given \ktree{k}{n - 1}\right] &= \ext{k}{n - 
		1} + \frac{k - 1}{(k - 1)(n - 1) + 1}\ext{k}{n - 1} + k
	\\ &= \frac{(k - 1)n + 1}{(k - 1)(n - 1) + 1}\ext{k}{n - 1} + k,
\end{align*}
equivalent to (after taking another expectation both sides)
\begin{equation}
	\label{Eq:extmeanrec}
	\E\left[\ext{k}{n}\right] = \frac{(k - 1)n + 1}{(k - 1)(n - 1) + 
		1}\E\left[\ext{k}{n - 1}\right] + k.
\end{equation}
Noting the initial condition $\ext{k}{0} = 0$, we solve 
Equation~(\ref{Eq:extmeanrec}) recursively for 
$\E\left[\ext{k}{n}\right]$ to get 
\begin{equation}
	\label{Eq:extmean}
	\E\left[\ext{k}{n}\right] = \frac{(kn - n + 1)m}{k - 1} 
	\left[\Psi\left(n + \frac{k}{k - 1}\right) - 
	\Psi\left(\frac{k}{k - 1}\right)\right],
\end{equation}
where $\Psi(\cdot)$ represents the {\em digamma function}. We 
finally apply a known formula for difference equation of digamma 
functions~\cite[Eq.\ 3.231.5]{gradshteyn2007table} to obtain the 
result 
stated in the theorem.

For the second moment of $\ext{k}{n}$, we consider the following 
convolution variable:
$$\dsq{k}{n} = \sum_{i = 1}^{L_n^{(k)}} \left(D^{(k)}_{i, 
	n}\right)^2,$$
the sum of squared depths of external nodes in $\ktree{k}{n}$. We 
implement a strategy analogous for the calculation of $\ext{k}{n}$ 
to compute the expectation of $\dsq{k}{n}$, which, later on, is used 
in the computation of the second moment of $\ext{k}{n}$. For better 
readability of the article, we omit the details of the derivation of 
$\E \left[ \dsq{k}{n} \right]$, but just state the result:
\begin{align}
	\E \left[ \dsq{k}{n} \right] &= \frac{k \bigl((k - 1)n + 
		1\bigr)}{k - 1} \left(\sum_{j = 1}^{n - 1}\left(\frac{2k 
		\Psi 
		\left(\frac{(k - 1)j + 2k - 1}{k - 1}\right) - 
		2k\Psi\left(\frac{k}{k - 1}\right)}{(k - 1)j + k} 
		\right.\right. 
	\nonumber
	\\ &\qquad{}+ \left. \left. \frac{(k - 1)\left((k - 1)j - 
		m\right)}{\left((k - 1)j + k\right)^2}\right)\right).
	\label{Eq:dsqexp}
\end{align}
Conditional on $\ktree{k}{n - 1}$ and label $i$, we square the 
almost-sure relation in Equation~(\ref{Eq:extinc}) to get
$$
\left(\ext{k}{n}\right)^2 = \left(\ext{k}{n - 1}\right)^2 + (k - 
1)^2 \left(D^{(k)}_{i, n - 1}\right)^2 + k^2 + 2 \ext{k}{n - 1} 
D^{(k)}_{i, n - 1} + 2k \ext{k}{n - 1} + 2k D^{(k)}_{i, n - 1}.
$$
We average out $i$ to obtain
\begin{align*}
	\E\left[\left(\ext{k}{n}\right)^2 \, \big{|} \, \ktree{n - 
		1}{k}\right] &= \left(1 + \frac{2(k - 1)}{(k - 1)(n - 1) + 
		1}\right) \left(\ext{k}{n - 1}\right)^2 + (k - 1)^2 
		\dsq{k}{n - 
		1}
	\\ &\qquad{} + \left(2k + \frac{2k}{(k - 1)(n - 1) + 1}\right) 
	\ext{k}{n - 1}.
\end{align*}
Taking the expectation with respect to $\ktree{k}{n - 1}$ and 
plugging in the results of $\E \left[ \dsq{k}{n} \right]$ (c.f.\ 
Equation~(\ref{Eq:dsqexp})) and $\E\left[\ext{k}{n}\right]$ (c.f.\ 
Equation~(\ref{Eq:extmean})), we obtain a recurrence for $\E 
\left[\left(\ext{k}{n}\right)^2 \right]$:
\begin{equation}
	\label{Eq:extsqrec}
	\E \left[\left(\ext{k}{n}\right)^2 \right] = \left(1 + \frac{2 
		(k - 1)}{(k - 1)(n - 1) + 1}\right) \E \left[\left(\ext{k}{n 
		- 
		1}\right)^2 \right]+ C(k, n),
\end{equation}
where $C(k, n)$ is a known function of $m$ and $n$. Solving the 
recurrence relation for $\E \left[\left(\ext{k}{n}\right)^2 \right]$ 
with the initial condition $\E \left[\left(\ext{k}{1}\right)^2 
\right] = m^2$, we obtain the stated result in the theorem.

At last, we present the exact expressions of the two functions: 
$E(k, n)$ and $C(k, n)$.
\begin{align*}
	E(k, n) &= \sum_{i = 1}^{n - 1} \frac{1}{\bigl((k - 1)i + 2k - 
		1\bigr)\bigl((k - 1)i + k\bigr)} \Biggl{\{}(2k^2 - 2k) 
		\sum_{j = 
		1}^{n - 1} \frac{\Psi \left(\frac{(k - 1)j + 2k - 1}{k - 
			1}\right)}{(k - 1)j + k} 
	\\ &\qquad{} - \left[2 \Psi\left(\frac{k}{k - 1}\right) - k + 
	1\right]\Psi \left(\frac{(k - 1)i + 2k - 1}{k - 1}\right)
	\\ &\qquad{} + 2k \left[ i \Psi\left(\frac{(k - 1)i + k}{k - 
		1}\right) + \Psi \left(1, \frac{(k - 1)i + 2k - 1}{k - 
		1}\right)\right]\Biggr{\}},
\end{align*}
where $\Psi(1, \cdot)$ is the first order derivative of the digamma 
function;
\begin{align*}
	C(k, n) &= \frac{2k}{k - 1} \Bigg{\{} k(k - 1)^2 \sum_{i = 0}^{n 
		- 1} \frac{\Psi \left(\frac{(k - 1)i + 2k - 1}{k - 
		1}\right)}{(k 
		- 1)i + k} - \left[(k^2 - k)n + \frac{k^2 + 
		1}{2}\right]\Psi\left(\frac{k}{k - 1}\right)
	\\ &\quad{} + (k^2 - k) \left[\Psi\left(1, \frac{(k - 1)n + 1}{k 
		- 1}\right) - \Psi\left(1, \frac{k}{k - 1}\right) + 
	\Psi^2\left(1, \frac{k}{k - 1}\right)\right]
	\\ &\quad{} + \Psi\left(\frac{(k - 1)n + 1}{k - 1}\right) 
	\left[(k^2 - k)n + \frac{k^2 + 1}{2} - (k^2 - 
	k)\Psi\left(\frac{k}{k - 1}\right)\right] + \frac{k^2 - k}{2} 
	\Bigg{\}}.
\end{align*}
\end{document}